\documentclass[reqno,11pt]{amsart}
\usepackage{amsmath}
\usepackage{amsxtra}
\usepackage{amscd}
\usepackage{amsthm}
\usepackage{amsfonts}
\usepackage{amssymb}
\usepackage{eucal}
\usepackage[matrix,arrow,curve]{xy}
\usepackage{tikz-cd}

\theoremstyle{plain}
\newtheorem{Thm}[subsection]{Theorem}
\newtheorem{Cor}[subsection]{Corollary}
\newtheorem{Lem}[subsection]{Lemma}
\newtheorem{Prop}[subsection]{Proposition}
\newtheorem{Conj}[subsection]{Conjecture}

\theoremstyle{definition}
\newtheorem{Def}[subsection]{Definition}

\theoremstyle{remark}

\newtheorem{Rem}[subsection]{Remark}

\newtheorem{conj}{Conjecture}

\numberwithin{equation}{section}
\renewcommand{\rm}{\normalshape}

\newif\ifShowLabels
\ShowLabelstrue
\newdimen\theight
\def\TeXref#1{%
    \leavevmode\vadjust{\setbox0=\hbox{{\tt
        \quad\quad  {\small \rm #1}}}%
    \theight=\ht0
    \advance\theight by \lineskip
    \kern -\theight \vbox to
    \theight{\rightline{\rlap{\box0}}%
    \vss}%
    }}%

\ShowLabelsfalse

\renewcommand{\sec}[2]{\section{#2}\label{S:#1}%
    \ifShowLabels \TeXref{{S:#1}} \fi}
\newcommand{\ssec}[2]{\subsection{#2}\label{SS:#1}%
    \ifShowLabels \TeXref{{SS:#1}} \fi}

\newenvironment{thm}[1]%
    { \begin{Thm} \label{T:#1}  \ifShowLabels \TeXref{T:#1} \fi }%
    { \end{Thm} }

\renewcommand{\th}[1]{\begin{thm}{#1} \sl }
\renewcommand{\eth}{\end{thm} }

\newenvironment{lemma}[1]%
    { \begin{Lem} \label{L:#1}  \ifShowLabels \TeXref{L:#1} \fi }%
    { \end{Lem} }
\newcommand{\lem}[1]{\begin{lemma}{#1} \sl}
\newcommand{\elem}{\end{lemma}}

\newenvironment{propos}[1]%
    { \begin{Prop} \label{P:#1}  \ifShowLabels \TeXref{P:#1} \fi }%
    { \end{Prop} }
\newcommand{\prop}[1]{\begin{propos}{#1}\sl }
\newcommand{\eprop}{\end{propos}}

\newenvironment{corol}[1]%
    { \begin{Cor} \label{C:#1}  \ifShowLabels \TeXref{C:#1} \fi }%
    { \end{Cor} }
\newcommand{\cor}[1]{\begin{corol}{#1} \sl }
\newcommand{\ecor}{\end{corol}}

\newenvironment{defeni}[1]%
    { \begin{Def} \label{D:#1}  \ifShowLabels \TeXref{D:#1} \fi }%
    { \end{Def} }
\newcommand{\defe}[1]{\begin{defeni}{#1} \sl }
\newcommand{\edefe}{\end{defeni}}

\newenvironment{remark}[1]%
    { \begin{Rem} \label{R:#1}  \ifShowLabels \TeXref{R:#1} \fi }%
    { \end{Rem} }
\newcommand{\rem}[1]{\begin{remark}{#1}}
\newcommand{\erem}{\end{remark}}

\newenvironment{conjec}[1]%
    { \begin{Conj} \label{Co:#1}  \ifShowLabels \TeXref{Co:#1} \fi }%
    { \end{Conj} }
\renewcommand{\conj}[1]{\begin{conjec}{#1} \sl }
\newcommand{\econj}{\end{conjec}}

\newcommand{\eq}[1]%
    { \ifShowLabels \TeXref{E:#1} \fi
       \begin{equation} \label{E:#1} }
\newcommand{\eeq}{ \end{equation} }

\newcommand{\prf}{ \begin{proof} }
\newcommand{\epr}{ \end{proof} }

\newcommand{\prft}{ \begin{proof} }
\newcommand{\eprt}{ \end{proof} }

\newcommand\nc{\newcommand}

\nc{\HC}{{\mathcal{HC}}}
\nc{\on}{\operatorname}
\nc{\BA}{{\mathbb{A}}}
\nc{\BC}{{\mathbb{C}}}
\nc{\BF}{{\mathbb{F}}}
\nc{\BG}{{\mathbb{G}}}
\nc{\BM}{{\mathbb{M}}}
\nc{\BN}{{\mathbb{N}}}
\nc{\BO}{{\mathbb{O}}}
\nc{\BQ}{{\mathbb{Q}}}
\nc{\BP}{{\mathbb{P}}}
\nc{\BR}{{\mathbb{R}}}
\nc{\BZ}{{\mathbb{Z}}}
\nc{\BS}{{\mathbb{S}}}

\nc{\CA}{{\mathcal{A}}}
\nc{\CB}{{\mathcal{B}}}
\nc{\CalC}{{\mathcal C}}
\nc{\CalD}{{\mathcal D}}
\nc{\CE}{{\mathcal{E}}}
\nc{\CF}{{\mathcal{F}}}
\nc{\CG}{{\mathcal{G}}}
\nc{\CH}{{\mathcal{H}}}
\nc{\CK}{{\mathcal{K}}}
\nc{\CL}{{\mathcal{L}}}
\nc{\CM}{{\mathcal{M}}}
\nc{\CMM}{{\mathcal{M}^{\operatorname{gen}}_\hbar(-\rho)}}
\nc{\CN}{{\mathcal{N}}}
\nc{\CO}{{\mathcal{O}}}
\nc{\CP}{{\mathcal{P}}}
\nc{\CQ}{{\mathcal{Q}}}
\nc{\CR}{{\mathcal{R}}}
\nc{\CS}{{\mathcal{S}}}
\nc{\CT}{{\mathcal{T}}}
\nc{\CU}{{\mathcal{U}}}
\nc{\CV}{{\mathcal{V}}}
\nc{\CW}{{\mathcal{W}}}
\nc{\CX}{{\mathcal{X}}}
\nc{\CY}{{\mathcal{Y}}}
\nc{\CZ}{{\mathcal{Z}}}

\nc{\gen}{{\operatorname{gen}}}
\nc{\cM}{{\check{\mathcal M}}{}}
\nc{\csM}{{\check{\mathcal A}}{}}
\nc{\obM}{{\overset{\circ}{\mathbf M}}{}}
\nc{\oCA}{{\overset{\circ}{\mathcal A}}{}}
\nc{\obA}{{\overset{\circ}{\mathbf A}}{}}
\nc{\ooM}{{\overset{\circ}{M}}{}}
\nc{\osM}{{\overset{\circ}{\mathsf M}}{}}
\nc{\vM}{{\overset{\bullet}{\mathcal M}}{}}
\nc{\nM}{{\underset{\bullet}{\mathcal M}}{}}
\nc{\obD}{{\overset{\circ}{\mathbf D}}{}}
\nc{\cp}{{\overset{\circ}{\mathbf p}}{}}
\nc{\ofZ}{{\overset{\circ}{\mathfrak Z}}{}}

\nc{\fa}{{\mathfrak{a}}}
\nc{\fb}{{\mathfrak{b}}}
\nc{\fg}{{\mathfrak{g}}}
\nc{\fgl}{{\mathfrak{gl}}}
\nc{\fh}{{\mathfrak{h}}}
\nc{\fj}{{\mathfrak{j}}}
\nc{\fm}{{\mathfrak{m}}}
\nc{\fn}{{\mathfrak{n}}}
\nc{\fu}{{\mathfrak{u}}}
\nc{\fp}{{\mathfrak{p}}}
\nc{\frr}{{\mathfrak{r}}}
\nc{\fs}{{\mathfrak{s}}}
\nc{\ft}{{\mathfrak{t}}}
\nc{\fT}{{\mathfrak{T}}}
\nc{\ofT}{{\overline{\mathfrak T}}}
\nc{\ofS}{{\overline{\mathfrak S}}}
\nc{\fsl}{{\mathfrak{sl}}}
\nc{\hsl}{{\widehat{\mathfrak{sl}}}}
\nc{\hgl}{{\widehat{\mathfrak{gl}}}}
\nc{\hg}{{\widehat{\mathfrak{g}}}}
\nc{\chg}{{\widehat{\mathfrak{g}}}{}^\vee}
\nc{\hn}{{\widehat{\mathfrak{n}}}}
\nc{\chn}{{\widehat{\mathfrak{n}}}{}^\vee}

\nc{\fA}{{\mathfrak{A}}}
\nc{\fB}{{\mathfrak{B}}}
\nc{\fD}{{\mathfrak{D}}}
\nc{\fE}{{\mathfrak{E}}}
\nc{\fF}{{\mathfrak{F}}}
\nc{\fG}{{\mathfrak{G}}}
\nc{\fI}{{\mathfrak{I}}}
\nc{\fJ}{{\mathfrak{J}}}
\nc{\fK}{{\mathfrak{K}}}
\nc{\fL}{{\mathfrak{L}}}
\nc{\fM}{{\mathfrak{M}}}
\nc{\fN}{{\mathfrak{N}}}
\nc{\frP}{{\mathfrak{P}}}
\nc{\fQ}{{\mathfrak Q}}
\nc{\fS}{{\mathfrak S}}
\nc{\fU}{{\mathfrak{U}}}
\nc{\fZ}{{\mathfrak{Z}}}

\nc{\ba}{{\mathbf{a}}}
\nc{\bb}{{\mathbf{b}}}
\nc{\bc}{{\mathbf{c}}}
\nc{\bd}{{\mathbf{d}}}
\nc{\be}{{\mathbf{e}}}
\nc{\bi}{{\mathbf{i}}}
\nc{\bj}{{\mathbf{j}}}
\nc{\bn}{{\mathbf{n}}}
\nc{\bp}{{\mathbf{p}}}
\nc{\br}{{\mathbf{r}}}
\nc{\bq}{{\mathbf{q}}}
\nc{\bu}{{\mathbf{u}}}
\nc{\bv}{{\mathbf{v}}}
\nc{\bx}{{\mathbf{x}}}
\nc{\by}{{\mathbf{y}}}
\nc{\bw}{{\mathbf{w}}}
\nc{\bA}{{\mathbf{A}}}
\nc{\bB}{{\mathbf{B}}}
\nc{\bC}{{\mathbf{C}}}
\nc{\bD}{{\mathbf{D}}}
\nc{\bE}{{\mathbf{E}}}
\nc{\bK}{{\mathbf{K}}}
\nc{\bH}{{\mathbf{H}}}
\nc{\bM}{{\mathbf{M}}}
\nc{\bN}{{\mathbf{N}}}
\nc{\bO}{{\mathbf{O}}}
\nc{\bQ}{{\mathbf Q}}
\nc{\bS}{{\mathbf{S}}}
\nc{\bT}{{\mathbf{T}}}
\nc{\bV}{{\mathbf{V}}}
\nc{\bW}{{\mathbf{W}}}
\nc{\bX}{{\mathbf{X}}}
\nc{\bP}{{\mathbf{P}}}
\nc{\bZ}{{\mathbf{Z}}}

\nc{\sA}{{\mathsf{A}}}
\nc{\sB}{{\mathsf{B}}}
\nc{\sC}{{\mathsf{C}}}
\nc{\sD}{{\mathsf{D}}}
\nc{\sF}{{\mathsf{F}}}
\nc{\sK}{{\mathsf{K}}}
\nc{\sM}{{\mathsf{M}}}
\nc{\sO}{{\mathsf{O}}}
\nc{\sQ}{{\mathsf{Q}}}
\nc{\sP}{{\mathsf{P}}}
\nc{\sV}{{\mathsf{V}}}
\nc{\sW}{{\mathsf{W}}}
\nc{\sZ}{{\mathsf{Z}}}
\nc{\sfp}{{\mathsf{p}}}
\nc{\sr}{{\mathsf{r}}}
\nc{\st}{{\mathsf{t}}}
\nc{\sfb}{{\mathsf{b}}}
\nc{\sfc}{{\mathsf{c}}}
\nc{\sd}{{\mathsf{d}}}
\nc{\sg}{{\mathsf{g}}}
\nc{\sk}{{\mathsf{k}}}
\nc{\sfl}{{\mathsf{l}}}

\nc{\BK}{{\bar{K}}}

\nc{\tA}{{\widetilde{\mathbf{A}}}}
\nc{\tB}{{\widetilde{\mathcal{B}}}}
\nc{\tg}{{\widetilde{\mathfrak{g}}}}
\nc{\tG}{{\widetilde{G}}}
\nc{\TM}{{\widetilde{\mathbb{M}}}{}}
\nc{\tN}{{\widetilde{\mathcal{N}}}{}}
\nc{\tO}{{\widetilde{\mathsf{O}}}{}}
\nc{\tU}{{\widetilde{\mathfrak{U}}}{}}
\nc{\TZ}{{\tilde{Z}}}
\nc{\tZ}{\widetilde{Z}{}}
\nc{\tx}{{\tilde{x}}}
\nc{\tbv}{{\tilde{\bv}}}
\nc{\tfP}{{\widetilde{\mathfrak{P}}}{}}
\nc{\tz}{{\tilde{\zeta}}}
\nc{\tmu}{{\tilde{\mu}}}

\nc{\td}{\ddot{\underline{d}}{}}
\nc{\tzeta}{\widetilde{\zeta}{}}
\nc{\hd}{{\widehat{\underline{d}}}}
\nc{\hG}{{\widehat{G}}}
\nc{\hBP}{\widehat{\mathbb P}{}}
\nc{\hQ}{{\widehat{Q}}}
\nc{\hsM}{\widehat{\mathsf M}{}}
\nc{\hfM}{\widehat{\mathfrak M}{}}
\nc{\hCP}{\widehat{\mathcal P}{}}
\nc{\hCR}{\widehat{\mathcal R}{}}
\nc{\hCS}{{\widehat{\mathcal S}}}
\nc{\hfZ}{\widehat{\mathfrak Z}{}}
\nc{\hZ}{\widehat{Z}{}}

\nc{\urho}{\underline{\rho}}
\nc{\uB}{\underline{B}}
\nc{\uC}{{\underline{\mathbb{C}}}}
\nc{\ui}{\underline{i}}
\nc{\ofP}{{\overline{\mathfrak{P}}}}

\nc{\hrho}{{\hat{\rho}}}

\nc{\unl}{\underline}
\nc{\ol}{\overline}
\nc{\one}{{\mathbf{1}}}
\nc{\two}{{\mathbf{t}}}

\nc{\Sym}{{\mathop{\operatorname{Sym}}}}
\nc{\Tot}{{\mathop{\operatorname{\normalshape Tot}}}}
\nc{\Hilb}{{\mathop{\operatorname{\normalshape Hilb}}}}
\nc{\Hom}{{\mathop{\operatorname{Hom}}}}
\nc{\CHom}{{\mathop{\operatorname{{\mathcal{H}}\it om}}}}
\nc{\defi}{{\mathop{\operatorname{\normalshape def}}}}
\nc{\length}{{\mathop{\operatorname{\normalshape length}}}}

\nc{\Cliff}{{\mathsf{Cliff}}}
\nc{\Fl}{{\mathcal{F}\ell}}
\nc{\Fib}{{\mathsf{Fib}}}
\nc{\Coh}{{\mathsf{Coh}}}
\nc{\FCoh}{{\mathsf{FCoh}}}

\nc{\reg}{{\text{\normalshape reg}}}
\nc{\res}{{\operatorname{res}}}
\nc{\cplus}{{\mathbf{C}_+}}
\nc{\cminus}{{\mathbf{C}_-}}
\nc{\cthree}{{\mathbf{C}_*}}
\nc{\Qbar}{{\bar{Q}}}

\nc{\bh}{{\bar{h}}}
\nc{\bOmega}{{\overline{\Omega}}}
\nc\tGr{\widetilde{\Gr}}

\nc{\seq}[1]{\stackrel{#1}{\sim}}
\nc\ogu{\overline{G/U}}
\nc\chlam{\check{\lam}}

\nc\St{\operatorname{St}}
\nc{\oZ}{{\overset{\circ}{Z}}}
\nc{\tF}{\widetilde{\mathcal F}}

\nc\uS{\underline{S}}
\nc\QM{\mathcal{QM}}

\nc{\chmu}{\check{\mu}}
\nc{\sm}{\setminus}
\nc{\ra}{\rightarrow}
\nc{\lar}{\leftarrow}
\nc{\La}{\Lambda}
\nc{\Lap}{\Lambda^{+}}
\nc{\oZal}{\overset{\circ}{Z^{\alpha}}}
\nc{\sig}{\sigma}
\nc{\al}{\alpha}
\nc{\la}{\lambda}
\nc{\is}{\simeq}
\nc{\ip}{\iota^{+}_{\la, \mu}}
\nc{\im}{\iota^{-}_{\la, \mu}}
\nc{\jp}{j^{+}_{\la, \mu}}
\nc{\jm}{j^{-}_{\la, \mu}}
\nc{\pip}{\pi^{+}_{\la, \mu}}
\nc{\pim}{\pi^{-}_{\la, \mu}}
\nc{\s}{\star}
\begin{document}

\author{Vasily Krylov}
\address{National Research University Higher School of Economics,
Russian Federation,
\newline Department of Mathematics, 6 Usacheva st, Moscow 119048}
\title[Integrable crystals and restriction to Levi via generalized slices]
{Integrable crystals and restriction to Levi via generalized slices
in the affine Grassmannian}
\email{kr-vas57@yandex.ru}

\begin{abstract}
Let $G$ be a connected reductive algebraic group over $\BC$.
Let $\Lap_{G}$ be the monoid of dominant weights of $G$.
We construct the integrable crystals $\bB^{G}(\la),\ \lambda\in\Lambda^+_G$,
using the geometry of generalized transversal slices in the affine Grassmannian
of the Langlands dual group. We construct the tensor product maps
$\bp_{\la_{1},\la_{2}}\colon \bB^{G}(\la_{1}) \otimes \bB^{G}(\la_{2}) \ra
\bB^{G}(\la_{1}+\la_{2})\cup\{0\}$  in terms of multiplication of generalized
transversal slices. Let $L \subset G$ be a Levi subgroup of $G$.
We describe the restriction to Levi
$\on{Res}^G_L\colon\on{Rep}(G)\to\on{Rep}(L)$ in terms of the hyperbolic
localization functors for the generalized transversal slices.
\end{abstract}
\maketitle

\sec{int}{Introduction}
Let us start with recalling the main results of Mirkovi\'c-Vilonen~\cite{MV}
and Beilinson-Drinfeld~\cite{BD} on geometric Satake isomorphism.

\ssec{Sat}{Geometric Satake isomorphism} Let us denote $\CK := \BC((z)), \CO := \BC[[z]]$, $\on{Gr}_{G}:=G(\CK)/G(\CO)$.
The tensor category $\on{Perv}_{G(\CO)}(\on{Gr}_{G})$ of $G(\CO)$-equivariant
perverse sheaves on $\on{Gr}_{G}$ was studied in~\cite{MV}. Let $\check{G}$ be the Langlands dual group. From \cite{MV} it follows that the categories $\on{Rep}(\check{G})$ and $\on{Perv}_{G(\CO)}(\on{Gr}_{G})$ are equivalent as tensor categories (let us denote the equivalence by $S_{G}: \on{Rep}(\check{G})\ra \on{Perv}_{G(\CO)}(\on{Gr}_{G})$). The corresponding fibre functor from $\on{Perv}_{G(\CO)}(\on{Gr}_{G})$ to $\on{Vect}_{\BC}$ sends a perverse sheaf to its global cohomology. We fix a Borel subgroup $B \subset G$ and a maximal torus $T \subset B$. Let $B_{-}$ denote the opposite Borel subgroup and $U$, $U_{-}$ be the unipotent radicals of $B$ and $B_{-}$. Let $\La_{G}$ be the coweight lattice of $T \subset G$ and let $\La_{G}^{+}$ be the submonoid of dominant coweights inside $\La_{G}$. Every representation $V$ of $\check{G}$ is graded by $\La_{G}$.
Accordingly, the fibre functor $\on{H}^{*}$ decomposes as the sum of  functors $\on{F}_{\mu}$, $\mu \in \La_{G}$. Let us describe the functors $\on{F}_{\mu}$.

One can identify $\La_{G}$ with $T(\CK)/T(\CO)$. Let us fix $\mu \in \La_{G}$ and let $z^{\mu}$ be the corresponding point inside $\on{Gr}_{G}$. The torus $T$ acts naturally on $\on{Gr}_G$ by the left multiplication. Let us consider the regular dominant coweight $2{\rho}_{G} : \BC^{*} \ra T$: the sum of positive coroots. It induces the $\BC^{*}$-action on $\on{Gr}_{G}$. Let us denote ${S}_{\mu}:=\{x \in \on{Gr}_{G} |  \on{lim}\limits_{t \ra 0} 2{\rho_{G}}(t)x=z^{\mu} \}=U(\CK) \cdot z^{\mu}$, ${T}_{\mu}:=\{x \in \on{Gr}_{G} |  \on{lim}\limits_{t \ra \infty} 2{\rho_{G}}(t)x=z^{\mu} \}=U_{-}(\CK) \cdot z^{\mu}.$

Let us consider the regular dominant weight $2\check{\rho}_{G} : T \ra \BC^{*}$: the sum of positive roots. The functor $\on{F}_{\mu}$ sends $\CF \in \on{Perv}_{G(\CO)}(\on{Gr}_{G})$ to the vector space $\on{H}^{2\check{\rho}_{G}(\mu)}_{{T}_{\mu}}(\on{Gr}_{G}, \CF)$. One can construct an isomorphism $\on{H}^{*}\simeq \bigoplus\limits_{\mu \in \La_{G}}\on{F}_{\mu}.$

\ssec{Wei}{Weights of irreducible representations} The irreducible objects of the category $\on{Perv}_{G(\CO)}(\on{Gr}_{G})$ are perverse sheaves $\on{IC}(\la):=j^{\la}_{!*}(\BC_{\on{Gr}^{\la}_{G}}[\on{dim}\on{Gr}^{\la}_{G}])$; here $\la \in \Lap_{G}$, $\on{Gr}^{\la}_{G}:=G(\CO)\cdot z^{\la}$ and $j^{\la}:\on{Gr}^{\la}_{G}\ra\on{Gr}_{G}$ denotes the locally closed embedding. Perverse sheaves $\on{IC}(\la)$ are supported on $\ol{\on{Gr}}^{\la}_G = \coprod \limits_{\mu \leq \la, \mu \in \Lap_{G}} \on{Gr}^{\mu}_G$. Let $V^{\la}_{\check{G}}$ denote the irreducible representation of $\check{G}$ of highest weight $\la$. The equivalence $S_{G}$ sends $V^{\la}_{\check{G}}$ to $\on{IC}(\la)$, so the $\check{T}$-weight component $(V^{\la}_{\check{G}})_{\mu}$ is isomorphic to $\on{F}_{\mu}(\on{IC}(\la))=\on{H}^{2\check{\rho}_{G}(\mu)}_{{T}_{\mu}}(\on{Gr}_{G}, \on{IC}(\la)).$ In \cite{MV} a canonical isomorphism $\on{H}^{2\check{\rho}_{G}(\mu)}_{{T}_{\mu}}(\on{Gr}_{G}, \on{IC}( \la)) \simeq \BC[\on{Irr}(T_{\mu} \cap \on{Gr}^{\la}_{G})]$ is constructed. Thus we get a distinguished basis in $(V^{\la}_{\check{G}})_{\mu}$ parametrized by irreducible components of $T_{\mu} \cap \on{Gr}^{\la}_{G}$ of the maximal dimension. Therefore the representation $V^{\la}_{\check{G}}$ has a basis parametrized by $\coprod \limits_{\mu \in \La_{G}} \on{Irr}(T_{\mu} \cap \on{Gr}^{\la}_{G})$.

\prop{Wei} \label{28} \cite{BG} Theorem 3.1.

(1) The set $\coprod \limits_{\mu \in \La_{G}} \on{Irr}(T_{\mu} \cap \on{Gr}^{\la}_{G})$ has a structure of $\check{G}$-crystal of highest weight $\la$.

(2) The crystals $\{ \coprod \limits_{\mu \in \La_{G}} \on{Irr}(T_{\mu} \cap \on{Gr}^{\la}_{G}) \}$ form a closed family.
\eprop

\ssec{TS}{Transversal slices}
Let us suppose $\mu \in \Lap_{G}$. Let us return to $ \on{H}^{2\check{\rho}_{G}(\mu)}_{{T}_{\mu}}(\on{Gr}_{G}, \on{IC}(\la))$. It can be identified with the hyperbolic restriction with respect to the $\BC^{*}$-action by $2{\rho}_{G}$ (see \cite{Br} for the definition of the hyperbolic restriction) of the sheaf $\on{IC}(\la)$ to the point $z^{\mu} \in \on{Gr}^{\mu}_{G}$. In \cite{KaTa} Kashiwara and Tanisaki defined the transversal slices $\ol{\CW}^{\la}_{\mu}$ to $\on{Gr}^{\mu}_{G}$ inside $\ol{\on{Gr}}^{\la}_{G}$ at the point $z^{\mu}$ (the slices are affine and conical with respect to the $\BC^{*}$-action by $2{\rho}_{G}$). One can see that the hyperbolic stalks of $\on{IC}$-sheaves on $\on{Gr}_{G}$ at the point $z^{\mu}$ are identified with the hyperbolic stalks of the corresponding $\on{IC}$-sheaves on transversal slices at the point $z^{\mu}$.

Thus, for every dominant $\mu$ we can define a locally closed affine subvariety $\ol{\CW}^{\la}_{\mu}$ in $\on{Gr}_{G}$ such that the irreducible components of repellents (under $2\rho_{G}$-action) $\ol{R}^{\la}_{\mu} \subset \ol{\CW}^{\la}_{\mu}$ to the point $z^{\mu}$ give us a basis in the space $(V^{\la}_{\check{G}})_{\mu}$. It is natural to try to describe crystal of highest weight $\bB^{G}(\la)$ in terms of transversal slices $\ol{\CW}^{\la}_{\mu}$ for various $\mu$. The problem is that transversal slices $\ol{\CW}^{\la}_{\mu}$ inside $\on{Gr}_{G}$ are defined only for dominant $\mu$. One can solve this problem by considering generalized transversal slices (see \cite[2(ii)]{BFN}) which are defined for any $\mu \in \La_{G}$.

Thus, for every $\la \in \Lap_{G}$, $\mu \in \La_{G}$ one has a closed affine subvariety $\ol{R}^{\la}_{\mu}$ inside a generalized slice $\ol{\CW}^{\la}_{\mu}$ such that for dominant $\mu$ the number of irreducible components of $\ol{R}^{\la}_{\mu}$ of maximal dimension is equal to $\on{dim}((V^{\la}_{\check{G}})_{\mu})$. In Section ~\ref{22} we prove (see Theorem ~\ref{7}~(1) below) that the natural morphism $\on{\bp}^{\la}_{\mu}:\ol{\CW}^{\la}_{\mu}\ra \ol{\on{Gr}}^{\la}_{G}$ restricts to
an isomorphism between $\ol{R}^{\la}_{\mu}$ and $T_{\mu} \cap \ol{\on{Gr}}^{\la}_{G}$. We give two different proofs of this statement. The first one uses an explicit "matrix" realization of a slice $\ol{\CW}^{\la}_{\mu}$ inside $G((z^{-1}))$ (see \cite[2($\on{xi}$)]{BFN}). The second proof uses the description of varieties $\ol{R}^{\la}_{\mu}$ and $T_{\mu} \cap \ol{\on{Gr}}^{\la}_{G}$ as moduli spaces.

So due to Proposition ~\ref{28}, the set $\coprod \limits_{\mu \in \La_{G}}\on{Irr}(\ol{R}^{\la}_{\mu})$ has a structure of $\check{G}$-crystal of highest weight $\la$ and the crystals $\{ \coprod \limits_{\mu \in \La_{G}}\on{Irr}(\ol{R}^{\la}_{\mu}) \}$ form a closed family. In Section ~\ref{32} we define the crystal structure on the set $\coprod \limits_{\mu \in \La_{G}}\on{Irr}(\ol{R}^{\la}_{\mu})$ using only the geometry of generalized transversal slices (mimicking \cite[Proposition 3.1]{BG}). From the previous observations it follows that for every $\la_{1},\la_{2}\in \Lap_{G}$, one has the retraction map between crystals $\on{\bp}_{\la_{1},\la_{2}}: \coprod\limits_{\mu \in \La_{G}}\on{Irr}(\ol{R}^{\la_{1}}_{\mu})\times \coprod\limits_{\mu \in \La_{G}}\on{Irr}(\ol{R}^{\la_{2}}_{\mu}) \ra \coprod\limits_{\mu \in \La_{G}}\on{Irr}(\ol{R}^{\la_{1}+\la_{2}}_{\mu})$. In Section ~\ref{23} we prove (see Theorem ~\ref{7}~(2) below) that these maps are induced by multiplication morphisms $\ol{\bf{\kappa}}^{\la_{1}, \la_{2}}_{\mu_{1}, \mu_{2}}:\ol{R}^{\la_{1}}_{\mu_{1}} \times \ol{R}^{\la_{2}}_{\mu_{2}} \ra \ol{R}^{\la_{1}+\la_{2}}_{\mu_{1}+\mu_{2}}$ between repellents. The proof
goes as follows: using results of Braverman-Gaitsgory and isomorphism from Theorem ~\ref{7} (1) we show that the maps $\on{\bp}_{\la_{1},\la_{2}}$ are induced by the convolution maps $\ol{\on{\bf{m}}}^{\la_{1},\la_{2}}_{\mu_{1},\mu_{2}}$ between $\on{Irr}((T_{\mu_{1}}\cap \ol{\on{Gr}}^{\la_{1}}_{G})\s (T_{\mu_{2}}\cap \ol{\on{Gr}}^{\la_{2}}_{G}))$ and $\on{Irr}((T_{\mu_{1}+\mu_{2}}\cap \ol{\on{Gr}}^{\la_{1}+\la_{2}}_{G}))$. After that we identify $(T_{\mu_{1}}\cap \ol{\on{Gr}}^{\la_{1}}_{G})\s (T_{\mu_{2}}\cap \ol{\on{Gr}}^{\la_{2}}_{G})$ with $\ol{R}^{\la_{1}}_{\mu_{1}} \times \ol{R}^{\la_{2}}_{\mu_{2}}$ and prove that under this identification the maps $\ol{\on{\bf{m}}}^{\la_{1},\la_{2}}_{\mu_{1},\mu_{2}}$ coincide with $\ol{\bf{\kappa}}^{\la_{1}, \la_{2}}_{\mu_{1}, \mu_{2}}$.

\ssec{RF}{Restriction to Levi}
\label{1}
Recall the results of Beilinson-Drinfeld~\cite[Proposition~5.3.29]{BD} on the geometric realization of the restriction to Levi $\on{Res}^{\check{G}}_{\check{L}}.$

Let $P_{-} \subset G$ be a parabolic subgroup of $G$ containing $B_{-}$. Let $L$ be the corresponding Levi factor. The diagram $G \hookleftarrow P_{-} \twoheadrightarrow L$ gives rise to two natural morphisms:

\begin{equation}
\xymatrix{ \on{Gr}_{G} & \ar[l]_{\iota} \on{Gr}_{P_{-}} \ar[r]^{\pi} & \on{Gr}_{L}
}
\end{equation}

It is known that the connected components of $\on{Gr}_{L}$ are parametrized by the character lattice $\La_{G,P_{-}}$ of the center $Z(\check{L})$ of $\check{L}$, and $\pi$ induces a bijection between the connected components of $\on{Gr}_{P_{-}}$ and the connected components of $\on{Gr}_{L}$. Let $\alpha_{G,P}: \La_{G} \ra \La_{G,P_{-}}$ denote the natural surjection.

For a character $\theta \in \La_{G,P_{-}}$ let us denote by $\on{Gr}_{P_{-},\theta}, \on{Gr}_{L,\theta}$ the corresponding connected components. Let $\iota_{\theta}, \pi_{\theta}$ denote the restrictions of $\iota$ and $\pi$ to the connected component $\on{Gr}_{P_{-},\theta}$. Let $\tilde{\theta}$ denote a lifting of $\theta$ to $\La_{G}$. Let $\check{\rho}_{G,L}:=\check{\rho}_{G}-\check{\rho}_{L}.$

Following~\cite[Proposition~5.3.29]{BD}, we consider the functor
$(\pi_{\theta})_{*}(\iota_{\theta})^{!}[2\check{\rho}_{G,L}(\tilde{\theta})]\colon
\on{Perv}_{G(\CO)}(\on{Gr}_{G})\to\on{Perv}_{L(\CO)}(\on{Gr}_{L})$.
The functor $R^{G,L}:=\bigoplus\limits_{\theta \in \La_{G,P_{-}}}(\pi_{\theta})_{*}(\iota_{\theta})^{!}[2\check{\rho}_{G,L}(\tilde{\theta})]$ coincides with the restriction functor $\on{Res}^{\check{G}}_{\check{L}}$ via identifications $S_{G}, S_{L}$, i.e. the following diagram is commutative:

\begin{equation}
\xymatrix{ \on{Perv}_{G(\CO)}(\on{Gr}_{G})  \ar[rr]^{R^{G,L}} &&  \on{Perv}_{L(\CO)}(\on{Gr}_{L})  \\ \\
\on{Rep}(\check{G}) \ar[uu]^{S_{G}}_{\wr} \ar[rr]^{\on{Res}^{\check{G}}_{\check{L}}} && \on{Rep}(\check{L}) \ar[uu]^{S_{L}}_{\wr}
}
\end{equation}

It is known that the natural grading on the functor $R^{G,L}$ corresponds to the $\La_{G,P_{-}}$-grading on the category $\on{Rep}(\check{L})$.

Let us denote $2{\rho}_{G,L}:=2{\rho}_{G}-2{\rho}_{L}$. Let us note that $2{\rho}_{G,L}$ is a "subregular" dominant (i.e. the corresponding element of $\La_{G}$ lies inside $\Lap_{G}$ and $\langle 2{\rho}_{G,L}, \al_{i} \rangle > 0$ for any $i \notin I_{L}$) cocharacter of the torus $Z(L)$. In Section ~\ref{24} we prove that the restriction functor $\on{Res}^{\check{G}}_{\check{L}}$ can be described via hyperbolic restriction functors from generalized slices $\ol{\CW}^{\la}_{\mu}$ to the fixed points subvarieties $(\ol{\CW}^{\la}_{\mu})^{2{\rho}_{G,L}}$ under the $\BC^{*}$-action by $2{\rho}_{G,L}$ (see Theorem ~\ref{2} for precise statement). For the proof we relate hyperbolic restriction functors $(\pi_{\theta})_{*}(\iota_{\theta})^{!}[2\check{\rho}_{G,L}(\tilde{\theta})]$ with the ones for slices. Lemma ~\ref{6} is the main observation which allows us to relate these two functors (i.e. to prove that one of them "restricts" to the other). This lemma is a natural generalization of Theorem ~\ref{7} (1).

\ssec{double}{Applications}
The above discussion shows that our constructions of integrable crystals and
restriction functors via generalized slices in the affine Grassmannian
essentially boil down to the well known existing constructions. We solved
this exercise with a view towards the (conjectural) slices in the double
affine Grassmannian ~\cite{BF} where it might prove useful for the geometric
constructions of integrable $\check{G}_{\on{aff}}$-crystals and the action of
$\check\fg_{\on{aff}}$ on the hyperbolic stalks of IC sheaves on the slices.

\ssec{Not}{List of notations}

We give the list of some notations here:

$\leq, \leq_{L}$: dominance orders with respect to $G$ or Levi subgroup $L$.

$W_{G}, W_{L}$: the Weyl groups corresponding to $T \subset G$ and $T \subset L$.

$w_{0}, {w_{0}}^{L}$: the longest elements of $W_{G}$ and $W_{L}$ respectively.

$\tilde{w}_{0}, {\tilde{w}_{0}}^{L}$ - representatives of  $w_{0}, {w_{0}}^{L}$ inside $G, L$ respectively.

$\check{G}, \check{L}$: the Langlands dual groups. $\check{T}$: the dual torus.

$g^{\la}_{G}:\ol{\on{Gr}}^{\la}_{G} \hookrightarrow \on{Gr}_{G}$, $g^{\nu}_{L}:\ol{\on{Gr}}^{\nu}_{L} \hookrightarrow \on{Gr}_{L}$,
$h^{\nu, \tilde{\la}}_{G}:\ol{\on{Gr}}^{\nu}_{L} \hookrightarrow \ol{\on{Gr}}^{\tilde{\la}}_{L}$: the closed embeddings.

$\ol{\CW}^{\la}_{\mu}$, $\ol{\CW}^{\nu}_{\mu,L}$: generalized transversal slices inside $\on{Gr}_{G}$, $\on{Gr}_{L}$.

$j^{\la}_{\mu}: \CW^{\la}_{\mu} \hookrightarrow \ol{\CW}^{\la}_{\mu}$: the natural open embedding.

$\ol{\bf{\kappa}}^{\la_{1}, \la_{2}}_{\mu_{1}, \mu_{2}}:\ol{R}^{\la_{1}}_{\mu_{1}} \times \ol{R}^{\la_{2}}_{\mu_{2}} \ra \ol{R}^{\la_{1}+\la_{2}}_{\mu_{1}+\mu_{2}}$: the multiplication morphism between repellents.

$B_{L}:=B\cap L$, $B_{L,-}:=B_{-}\cap L$, $U_{P_{-}}$: the unipotent radical of $P_{-}$.

$\varkappa^{\la}_{\mu,G,L}:(\ol{\CW}^{\la}_{\mu})^{2{\rho}_{G,L}} \hookrightarrow \ol{\CW}^{\tilde{\la}}_{\mu,L}$, $\imath^{\nu}_{\mu,G,L}: \ol{\CW}^{\nu}_{\mu,L} \hookrightarrow (\ol{\CW}^{\la}_{\mu})^{2{\rho}_{G,L}}$,
$\ell^{'\la,\la}_{\mu}:\ol{\CW}^{'\la}_{\mu} \hookrightarrow \ol{\CW}^{\la}_{\mu}$, $\ell^{'\nu,\nu}_{\mu,L}:\ol{\CW}^{'\nu}_{\mu,L} \hookrightarrow \ol{\CW}^{\nu}_{\mu,L}$: the closed embeddings.

$\bp^{\la}_{\mu}:\ol{\CW}^{\la}_{\mu}\ra \ol{\on{Gr}}^{\la}_{G}$,
$\bp^{\tilde{\la}}_{\mu,L}:\ol{\CW}^{\tilde{\la}}_{\mu,L}\ra \ol{\on{Gr}}^{\tilde{\la}}_{L}$: the natural morphisms.

For a linear algebraic group $H$ we will denote by $H[z^{-1}]_{1}, H[[z^{-1}]]_{1}$ kernels of the natural ("evaluation at $\infty$") homomorphisms $H[z^{-1}] \ra H, H[[z^{-1}]]\ra H$.

$\on{\bf{mult}}: U_{-}[z^{-1}]_{1} \times U_{-}[[z]] \ra U_{-}((z))$: the multiplication morphism.

$\xi_{U_{P_{-}}}:$ $U_{P_{-}}[z^{-1}]_{1}\tilde{\ra} U_{P_{-}}(\CK)/U_{P_{-}}(\CO)$: the natural isomorphism.

$\on{Irr}(X)$: the set of irreducible components of the maximal dimension of a variety $X$.

\rem{RF}

For a scheme $X$ over $\BC$ consider the corresponding variety $X_{\on{red}}$. We will use the following observation: suppose the following square of schemes of finite type over $\BC$ is cartesian:

\begin{equation}
\xymatrix{ P  \ar[r] \ar[d] &  Y \ar[d] \\
X \ar[r]& S
}
\end{equation}
then the following square of varieties is cartesian {\em in the category of
complex varieties}:

\begin{equation}
\xymatrix{ P_{\on{red}}  \ar[r] \ar[d] &  Y_{\on{red}} \ar[d] \\
X_{\on{red}} \ar[r]& S_{\on{red}}
}
\end{equation}
\erem

\ssec{Int}{Organization of the paper} In Section ~\ref{41} we give the definitions of the main objects of our study. In Section ~\ref{42} we formulate the main theorems of the paper. Section ~\ref{22} is devoted to the proof of Theorem ~\ref{7} (1).  In Section ~\ref{24} we prove Theorem ~\ref{2}.  In section ~\ref{32} we recall the results and the main constructions of \cite{BG} and prove the first part of Theorem ~\ref{7} (2). Section ~\ref{23} is devoted to the proof of the second part of Theorem ~\ref{7} (2).

\ssec{Int}{Acknowledgements} I would like to thank Michael Finkelberg for posing
the problem and for many helpful discussions and numerous suggestions. The author was supported in part by Laboratory of Algebraic Geometry and its Applications, August M\"obius contest (2016) and Dobrushin stipend.

\sec{Def}{Definitions} \label{41}

\ssec{Gr}{Affine Grassmannian} Recall that $G$ is a connected reductive algebraic group over $\BC$. Let $\CK := \BC((z)), \CO := \BC[[z]]$. Let $\on{Gr}_G := G(\CK)/G(\CO)$. The affine Grassmannian $\on{Gr}_G$ is the set of $\BC$-points of an ind-scheme over $\BC$ which we will denote by the same symbol. One can think about $\on{Gr}_{G}$ as the moduli space of principal $G$-bundles $\CP$ on $\BP^{1}$ with a trivialization $\sig: \CP_{triv}|_{\BP^{1} \sm 0} \tilde{\ra} \CP|_{\BP^{1} \sm 0}$. We fix a maximal torus $T \subset G$. Recall that $\La_{G}$ is the coweight lattice of $T\in G$ and $\Lap_{G}$ is the submonoid of the dominant coweights. One can identify $\La_{G}$ with $T(\CK)/T(\CO)$. Fix $\la \in \La_{G}$ and let $z^{\la}$ denote any lift of $\la$ to $T(\CK)$. We will denote by the same symbol the corresponding point inside $\on{Gr}_{G}$. The group $G(\CO)$ acts on $\on{Gr}_G$ on the left. Let  $\on{Gr}^{\la}_{G}:=G(\CO) \cdot z^{\la}$ (with reduced scheme structure). It is known that $\on{Gr}_{G}=\coprod \limits_{\la \in \Lap_{G}} \on{Gr}^{\la}_G$. Also $\ol{\on{Gr}}^{\la}_G = \coprod \limits_{\mu \leq \la, \mu \in \Lap_{G}} \on{Gr}^{\mu}_{G}$.

\rem{Gr} It follows that for $\mu \in \La_{G}$, $\on{Gr}^{\mu}_G \subset \ol{\on{Gr}}^{\la}_{G}$  if and only if $\mu$ is a weight of the irreducible representation $V^{\la}_{\check{G}}$ of $\check{G}$.
\erem

Let $T^{\la}_{\mu}:=T_{\mu}\cap \on{Gr}^{\la}_{G}$, $\ol{T}^{\la}_{\mu}:=T_{\mu}\cap \ol{\on{Gr}}^{\la}_{G}.$

\ssec{Conv}{Convolution} Let us define the ind-scheme $\on{Gr}_{G}\s\on{Gr}_{G}$
as $G(\CK) \times_{G(\CO)} \on{Gr}_{G}$. Let $\pi: G(\CK) \times \on{Gr}_{G} \ra \on{Gr}_{G}\s\on{Gr}_{G}$ denote the natural projection. Let $p_{1}, p_{2}: G(\CK) \times \on{Gr}_{G} \ra \on{Gr}_{G}$ and $\on{\bf{m}}: \on{Gr}_{G}\s\on{Gr}_{G} \ra \on{Gr}_{G}$ be defined as follows. Take $g \in G(\CK), [x] \in \on{Gr}_{G}$; then $p_{1}(g,[x]) := [g]$, $p_{2}(g,[x]) := [x]$, $\on{\bf{m}}(g,[x])=[gx].$

Let $\on{Gr}^{\la_{1}}_{G}\s \on{Gr}^{\la_{2}}_{G}:=\pi(p_{1}^{-1}(\on{Gr}^{\la_{1}}_{G}) \cap p_{2}^{-1}(\on{Gr}^{\la_{2}}_{G})),$

$(\on{Gr}^{\la_{1}}_{G}\s\on{Gr}^{\la_{2}}_{G})^{\la_{3}}:=\on{\bf{m}}^{-1}(\on{Gr}^{\la_{3}}_{G})\cap (\on{Gr}^{\la_{1}}_{G}\s\on{Gr}^{\la_{2}}_{G}),$

$\ol{\on{Gr}}^{\la_{1}}_{G}\s \ol{\on{Gr}}^{\la_{2}}_{G}:=\pi(p_{1}^{-1}(\ol{\on{Gr}}^{\la_{1}}_{G}) \cap p_{2}^{-1}(\ol{\on{Gr}}^{\la_{2}}_{G})).$

\ssec{Zast}{Quasi-maps to flag variety and Zastava} We fix a Borel and a Cartan subgroup $G \supset B \supset T$. Let us denote by $I_{G}$ the set parameterizing simple coroots of $G$. Let us denote by $\La^{\on{pos}}_{G}$ the submonoid of $\La_{G}$ spanned by the simple coroots $\alpha_{i}, i \in I_{G}$. For $\al \in \La^{\on{pos}}_{G}$ let us denote by $\on{QMaps}^{\al}(\BP^{1},\CB)$ the moduli space of quasimaps $\phi$ of degree $\al$ from $\BP^{1}$ to the flag variety $\CB := G/B$ (see \cite[Subsection 2.2]{B}). Let $ \check{\La}_{G}^{+}$ be a monoid of dominant weights of $T \subset G$. One can think about a point in $\on{QMaps}^{\al}(\BP^{1},\CB)$ as a collection of invertible subsheaves $\CL_{\check{\eta}} \subset V^{\check{\eta}}_{G} \otimes \CO_{\BP^{1}}$ of degree $-\langle \al, \check{\eta} \rangle$ for each $\check{\eta} \in \check{\La}_{G}^{+}$ subject to Pl$\ddot{\on{u}}$cker relations.

For $\al \in \La^{\on{pos}}_{G}$ let us denote by $Z^{\al} \subset \on{QMaps}^{\al}(\BP^{1},\CB)$ the moduli space of quasimaps $\phi$ of degree $\al$ from $\BP^{1}$ to the flag variety $\CB := G/B$ such that $\phi$ has no defect at $\infty \in \BP^{1}$ and $\phi(\infty)=B_{-}$. It contains the open dense moduli subspace $\oZal$ of based maps.
\ssec{Gensl}{Generalized slices and repellents} \label{1234567} (see \cite[2(ii)]{BFN}) Let $\la$ be a dominant coweight of $G$ and let $\mu \leq \la$ be an arbitrary coweight of $G$. Let $\al := \la - \mu$. Let us denote by $\ol{\CW}^{\la}_{\mu}(\BC)$ the moduli space of the following data:

(a) G-bundle $\CP$ on $\BP^{1}.$

(b) A trivialization  $\sig: \CP_{triv}|_{\BP^{1} \sm 0} \tilde{\ra} \CP|_{\BP^{1} \sm 0}$ having a pole of degree $\leq \la$. This means that for an irreducible G-module $V^{\check{\la}}_{G}$ and the associated vector bundle $\CV_{\CP}^{\check{\la}}$ we have $V^{\check{\la}}_{G} \otimes \CO_{\BP^{1}}(-\langle \la, \check{\la} \rangle \cdot 0) \subset \CV_{\CP}^{\check{\la}} \subset V^{\check{\la}}_{G} \otimes \CO_{\BP^{1}}(-\langle w_{0}(\la), \check{\la} \rangle \cdot 0)$. In other words this means that the point $(\CP, \sig) \in \on{Gr}_G$ lies inside $\ol{\on{Gr}}_G^{\la}$.

(c) A $B$-structure $\phi$ on $\CP$ of degree $w_{0}(\mu)$ having no defect at $\infty$ and having fiber $B_{-}$ at $\infty$ (with respect to $\sig$). In Pl\"ucker coordinates it means that for every $\check{\eta} \in \check{\La}^{+}_{G}$ we have an invertible subsheaf $\CL_{\check{\eta}} \subset \CV_{\CP}^{\check{\eta}}$ of degree $-\langle w_{0}(\mu), \check{\eta} \rangle$.

Let us endow $\ol{\CW}^{\la}_{\mu}(\BC)$ with the variety structure in the following way: $$\ol{\CW}^{\la}_{\mu}:=(\ol{\on{Gr}}^{\la}_{G}\times_{'\on{Bun}_{G}(\BP^{1})}\on{Bun}^{w_{0}\mu}_{B}(\BP^{1}))_{\on{red}},$$
where $'\on{Bun}_{G}(\BP^{1})$ is the moduli stack of $G$-bundles on $\BP^{1}$ with a $B$-structure at $\infty$ and $\on{Bun}^{w_{0}\mu}_{B}(\BP^{1})$ is the moduli stack of $B$-bundles on $\BP^{1}$ of degree $w_{0}\mu$.

\rem{gensl} Following \cite{BFN} let us note that the scheme $\ol{\on{Gr}}^{\la}_{G}\times_{'\on{Bun}_{G}(\BP^{1})}\on{Bun}^{w_{0}\mu}_{B}(\BP^{1})$ is actually reduced since it is generically reduced and Cohen-Macaulay (see \cite[Lemma 2.16]{BFN}).
\erem

\rem{gensl} Actually variety $\ol{\CW}^{\la}_{\mu}$ can be defined for any $\la \in \Lap_{G}$, $\mu \in \La_{G}$ but it will be nonempty precisely if $\mu \leq \la$.
\erem

We have a convolution diagram
\begin{equation}
\xymatrix{
\ol{\on{Gr}}^{\la}_{G} &&  \ar[ll]_{\on{\bf{p}}^{\la}_{\mu}} \ol{\CW}^{\la}_{\mu}  \ar[rr]^{\on{\bf{q}}^{\la}_{\mu}} && Z^{-w_{0}(\la-\mu)},
}
\end{equation}
where $\on{\bf{p}}^{\la}_{\mu}$ maps $(\CP, \sig, \phi)$ to $(\CP, \sig)$ and $\on{\bf{q}}^{\la}_{\mu}$ maps $(\CP, \sig, \phi)$ to the collection of subsheaves $\CL_{\check{\eta}}(\langle w_{0}(\la), \check{\eta} \rangle \cdot 0) \subset V^{\check{\eta}} \otimes \CO_{\BP^{1}}$ where the embedding is induced by $\sig^{-1}$.

By definition, the morphism $\on{\bf{p}}^{\la}_{\mu} \times \on{\bf{q}}^{\la}_{\mu} : \ol{\CW}^{\la}_{\mu} \ra \ol{\on{Gr}}^{\la}_G \times Z^{-w_{0}(\la-\mu)}$ is a locally closed embedding.

Let us define ${\CW}^{\la}_{\mu}:=(\on{\bp}^{\la}_{\mu})^{-1}(\on{Gr}^{\la}_{G}).$ Let us denote by $j^{\la}_{\mu}:{\CW}^{\la}_{\mu} \ra \ol{\CW}^{\la}_{\mu}$ the corresponding open embedding. We can note that for every dominant coweight $'\la\leq\la$ one has the natural closed embedding $\ell^{'\la,\la}_{\mu}:\ol{\CW}^{'\la}_{\mu} \ra \ol{\CW}^{\la}_{\mu}$.

Let us describe the actions of the Cartan torus $T$ on the varieties $\ol{\CW}^{\la}_{\mu}$, $\on{Gr}^{\la}_{G}$, $Z^{-w_{0}(\la-\mu)}$ such that morphisms $\on{\bf{p}}^{\la}_{\mu}$, $\on{\bf{q}}^{\la}_{\mu}$ are $T$-equivariant with respect to these actions: $T$ acts on $\ol{\CW}^{\la}_{\mu}$, $\on{Gr}^{\la}_{G}$ by changing trivialization, $T$ acts on $Z^{-w_{0}(\la-\mu)}$ via its natural action on $\CB$. In particularly we get the $\BC^{*}$-action on $\ol{\CW}^{\la}_{\mu}$ via $2{\rho}_{G}: \BC^{*} \ra T$.

\lem{Gensl} \label{100001} The subvariety $(\ol{\CW}^{\la}_{\mu})^{2{\rho}_{G}} \subset \ol{\CW}^{\la}_{\mu}$ consists of one point if $\mu$ is a weight of $V^{\la}_{\check{G}}$ and is empty otherwise.
\elem

\prf It is enough to show it on the level of $\BC$-points. Now the statement directly follows from the "matrix" description of slices in Section ~\ref{21} below:

$\ol{\CW}^{\la}_{\mu}(\BC)^{2{\rho}_{G}}=(B[[z^{-1}]]_{1}z^{\mu}B_{-}[[z^{-1}]]_{1})^{2{\rho}_{G}}\bigcap \ol{G[z]z^{\la}G[z]}=$

$T[[z^{-1}]]_{1}z^{\mu}T[[z^{-1}]]_{1}\bigcap \ol{G[z]z^{\la}G[z]}=T[[z^{-1}]]_{1}z^{\mu}T[[z^{-1}]]_{1}\bigcap \coprod\limits_{\mu' \leq \la, \mu' \in \Lap_{G}}G[z]z^{\mu'}G[z]=$

$T[z^{-1}]_{1}z^{\mu}T[z^{-1}]_{1}\bigcap \coprod\limits_{\mu' \leq \la, \mu' \in \Lap_{G}}G[z]z^{\mu'}G[z]=z^{\mu}\bigcap \coprod\limits_{\mu' \leq \la, \mu' \in \Lap_{G}}G[z]z^{\mu'}G[z].$

Let us prove that $\ol{G[z]z^{\la}G[z]}=\coprod\limits_{\mu' \leq \la, \mu' \in \Lap_{G}}G[z]z^{\mu'}G[z].$ Let us fix an embedding $\imath:G \hookrightarrow GL_{m}$. There exists a big enough positive integer $N$ such that for any $g \in G[z]z^{\la}G[z]$, the coefficients of the matrix $z^{N}\imath(g)$ are polynomials in $z$. Thus the set $G[z]z^{\la}G[z]$ is contained in the subset consisting of all $g \in G((z^{-1}))$ such that coefficients of $z^{N}\imath(g)$ are polynomials in $z$
(we denote this subset $G((z^{-1}))_{N}$). The ind-scheme $G((z^{-1}))_{N}$ is obviously closed inside $G((z^{-1}))$, and $G[z]z^{\la}G[z] \subset G((z^{-1}))_{N}$,
so $\ol{G[z]z^{\la}G[z]} \subset G((z^{-1}))_{N}.$

Let us consider the natural projection $\on{\bf{pr}}:G(z)_{N}\ra \on{Gr}_{G}$. The image of $G[z]z^{\la}G[z]$ lies in the subset $\on{Gr}^{\la}_{G}$. The closure $\ol{\on{Gr}}^{\la}_{G}=\coprod\limits_{\mu' \leq \la, \mu' \in \Lap_{G}}\on{Gr}^{\mu'}_{G}$ thus $\coprod\limits_{\mu' \leq \la, \mu' \in \Lap_{G}}G[z]z^{\mu'}G[z]=\on{\bf{pr}}^{-1}(\coprod\limits_{\mu' \leq \la, \mu' \in \Lap_{G}}\on{Gr}^{\mu'}_{G})$ is closed (the last equality holds by \cite[Proposition 1.2.4]{Gi}).

So the closure $\ol{G[z]z^{\la}G[z]}$ lies in $\coprod\limits_{\mu' \leq \la, \mu' \in \Lap_{G}}G[z]z^{\mu'}G[z]$.

To prove the embedding in the opposite direction it is enough to show that for every dominant $\mu' \leq \la$, $z^{\mu'} \in \ol{G[z]z^{\la}G[z]}$. For dominant $\mu'$ the transversal slice $\ol{\CW}^{\la}_{\mu'}$ coincides with the corresponding transversal slice inside $\on{Gr}_{G}$, thus $z^{\mu'} \in \ol{\CW}^{\la}_{\mu'} \subset \ol{G[z]z^{\la}G[z]}$.
\epr

Let $\ol{R}^{\la}_{\mu} \subset \ol{\CW}^{\la}_{\mu}$ be the repellent to the point $z^{\mu}$ inside $\ol{\CW}^{\la}_{\mu}$ with respect to the $\BC^{*}$-action on $\ol{\CW}^{\la}_{\mu}$ via $2{\rho_{G}}$. In more details let us denote $\ol{R}^{\la}_{\mu}:=\{x \in \ol{\CW}^{\la}_{\mu} | \on{lim}\limits_{t \ra \infty} 2{\rho_{G}}(t)x=z^{\mu}\}$. The map $\on{\bf{p}}^{\la}_{\mu}$ is $T$-equivariant, so $\on{\bf{p}}^{\la}_{\mu}$ maps $\ol{R}^{\la}_{\mu}$ to $\ol{T}^{\la}_{\mu}$.

\ssec{Mult}{Multiplication of slices.} (see \cite[2($\on{vi}$)]{BFN}) Let $\la_{1}, \la_{2} \in \Lap_{G}$, $\mu_{1}, \mu_{2} \in \La_G$. Let us construct a multiplication morphism $\ol{\CW}^{\la_{1}}_{\mu_{1}} \times \ol{\CW}^{\la_{2}}_{\mu_{2}} \ra \ol{\CW}^{\la_{1}+\la_{2}}_{\mu_{1}+\mu_{2}}$. Let us describe a symmetric definition of generalized slices. Let us take $\mu_{-}, \mu_{+} \in \La_{G}$. Let us denote by $\ol{\CW}^{\la}_{\mu_{-}, \mu_{+}}$ the moduli space of the following data:

(a) $G$-bundles $\CP_{-}, \CP_{+}$ on $\BP^{1}$.

(b) An isomorphism $\sig: \CP_{-}|_{\BP^{1}\setminus 0} \tilde{\ra} \CP_{+}|_{\BP^{1}\setminus 0}$ having a pole of degree $\leq \la$ at $0 \in \BP^{1}$.

(c) A trivialization of $\CP_{-}=\CP_{+}$ at $\infty$.

(d) A $B_{-}$ structure $\phi_{-}$ on $\CP_{-}$ such that the induced $T$-bundle has degree $-w_{0}(\mu_{-})$ and the fiber of $\phi_{-}$ at $\infty$ is $B$.

(e) A $B_{+}$ structure $\phi_{+}$ on $\CP_{+}$ such that the induced $T$-bundle has degree $w_{0}(\mu_{+})$ and the fiber of $\phi_{+}$ at $\infty$ is $B_{-}$.

\prop{Mult} (see \cite[2($\on{v}$)]{BFN}) We have canonical isomorphisms
$\ol{\CW}^{\la}_{\mu_{-}, \mu_{+}} \simeq \ol{\CW}^{\la}_{\mu_{-}+\mu_{+}, 0} \simeq \ol{\CW}^{\la}_{\mu_{-}+\mu_{+}}$.
\eprop

Thus to define the multiplication morphism it is enough to construct a morphism

 $\ol{\CW}^{\la_{1}}_{\mu_{1}, 0} \times \ol{\CW}^{\la_{2}}_{0, \mu_{2}} \ra \ol{\CW}^{\la_{1}+\la_{2}}_{\mu_{1}, \mu_{2}}$.

It is given by $((\CP^{1}_{-},\CP^{1}_{+},\sig_{1},\phi^{1}_{-},\phi^{1}_{+}),(\CP^{2}_{-},\CP^{2}_{+},\sig_{2},\phi^{2}_{-},\phi^{2}_{+})) \mapsto (\CP^{1}_{-}, \CP^{2}_{+}, \sig_{2} \circ \sig_{1}, \phi^{1}_{-}, \phi^{2}_{+})$.

\ssec{Cry}{Crystals} (see \cite{Ka}) A $\check{\mathfrak{g}}$-crystal is a set $\bB$ together with maps:

  (1) $wt: \bB \ra \La_{G}, \hspace{0,1cm} \varepsilon_{i}, \hspace{0,1cm}  \varphi_{i}: \bB \ra \mathbb{Z} \cup \infty.$

  (2) $\on{e}_{i}, \hspace{0,1cm} \on{f}_{i}: \bB \ra \bB \cup \{ 0\}.$

such that for each $i \in I_{G}$ we have:

(a) For each $\bb \in \bB$, \hspace{0,1cm}  $\varphi_{i}(\bb) = \varepsilon_{i}(\bb) + \langle wt(\bb), \alpha_{i}^\vee \rangle.$

(b) Let $\bb \in \bB$. If $\on{e}_{i} \cdot \bb \in \bB$ for some $i$, then

$$wt(\on{e}_{i} \cdot \bb) = wt(\bb) + \alpha_{i}, \hspace{0,1cm}  \varepsilon_{i}(\on{e}_{i} \cdot \bb) = \varepsilon_{i}(\bb) - 1, \hspace{0,1cm}  \varphi_{i}(\on{e}_{i} \cdot \bb) = \varphi_{i}(\bb) + 1.$$

 If $\on{f}_{i} \cdot \bb \in \bB$ for some $i$, then

$$wt(\on{f}_{i} \cdot \bb) = wt(\bb) - \alpha_{i}, \hspace{0,1cm}  \varepsilon_{i}(\on{f}_{i} \cdot \bb) = \varepsilon_{i}(\bb) + 1, \hspace{0,1cm}  \varphi_{i}(\on{f}_{i} \cdot \bb) = \varphi_{i}(\bb) - 1.$$

(c) For all $\bb, \hat{\bb} \in \bB$, \hspace{0,1cm}  $\on{e}_{i} \cdot \bb = \hat{\bb}$ if and only if $\bb = \on{f}_{i} \cdot \hat{\bb}$.

A $\check{\mathfrak{g}}$-crystal is called normal if for any $\bb \in \bB, i \in I_{G}$, we have

$$\varepsilon_{i}(\bb) = \on{max} \{n \in \mathbb{N}| (\on{e}_{i})^{n}(\bb) \neq 0 \}, \hspace{0,1cm}  \varphi_{i}(\bb) = \on{max} \{n \in \mathbb{N}| (\on{f}_{i})^{n}(\bb) \neq 0 \}.$$

\ssec{CryT}{Tensor product of crystals} Given two $\check{\mathfrak{g}}$-crystals $\bB, \hat{\bB}$ one can define a $\check{\mathfrak{g}}$-crystal structure on the set $\bB \times \hat{\bB}$:

$$wt(\bb \otimes \hat{\bb}) = wt(\bb) + wt(\hat{\bb})$$

 \begin{equation} \on{e}_{i} \cdot (\bb \otimes \hat{\bb}) =
\begin{cases}
\on{e}_{i} \cdot \bb \otimes \hat{\bb}, \hspace{0,1cm} \on{if} \hspace{0,1cm}  \varepsilon_{i}(\bb) > \varphi_{i}(\hat{\bb})\\
\bb \otimes \on{e}_{i} \cdot \hat{\bb}, \hspace{0,1cm}  \on{otherwise}
\end{cases}
\end{equation}

 \begin{equation} \on{f}_{i} \cdot (\bb \otimes \hat{\bb}) =
\begin{cases}
\on{f}_{i} \cdot \bb \otimes \hat{\bb}, \hspace{0,1cm}  \on{if} \hspace{0,1cm}  \varepsilon_{i}(\bb) \ge \varphi_{i}(\hat{\bb})\\
\bb \otimes \on{f}_{i} \cdot \hat{\bb}, \hspace{0,1cm}  \on{otherwise}
\end{cases}
\end{equation}

$$\varepsilon_{i}(\bb \otimes \hat{\bb}) = \on{max} \{ \varepsilon_{i}(\hat{\bb}), \hspace{0,1cm}  \varepsilon_{i}(\bb) - \varphi_{i}(\hat{\bb}) + \varepsilon_{i}(\hat{\bb}) \}$$

$$\varphi_{i}(\bb \otimes \hat{\bb}) = \on{max} \{ \varphi_{i}(\bb), \hspace{0,1cm}  \varphi_{i}(\hat{\bb}) - \varepsilon_{i}(\bb) + \varphi_{i}(\bb) \}.$$

Let us denote this crystal $\bB \otimes \hat{\bB}$. It is known that if $\bB, \hat{\bB}$ are normal then $\bB \otimes \hat{\bB}$ is normal.

\ssec{MPr}{Morphisms of crystals}(see \cite[Section 7]{Ka}) Let $\bB_{1}, \bB_{2}$ be two $\check{\mathfrak{g}}$-crystals. Let $\bp: \bB_{1} \ra \bB_{2} \cup \{0\}$ be a map from $\bB_{1}$ to $\bB_{2} \cup \{0\}$. Let us say that $\bp$ is a morphism of $\check{\mathfrak{g}}$-crystals $\bB_{1}$ and $\bB_{2}$ if the following conditions hold:

(a) $\bp$ commutes with $wt, \varepsilon_{i}, \varphi_{i}$ for any $i \in I_{G}$,

(b) $\bp(\on{e}_{i} \cdot \bb) = \on{e}_{i} \cdot \bp(\bb)$ for any $\bb \in \bB_{1}, i \in I_{G}$ such that $\bp(\bb) \neq 0$ and $\bp(\on{e}_{i} \cdot \bb) \neq 0$,

(c) $\bp(\on{f}_{i} \cdot \bb) = \on{f}_{i} \cdot \bp(\bb)$ for any $\bb \in \bB_{1}, i \in I_{G}$ such that $\bp(\bb) \neq 0$ and $\bp(\on{f}_{i} \cdot \bb) \neq 0$.

\quad

Let us say that $\bp$ is a strict morphism of $\check{\mathfrak{g}}$-crystals $\bB_{1}$ and $\bB_{2}$ if $\bp$ is a morphism of crystals and $\bp$ commutes with $\on{e}_{i}, \on{f}_{i}$ for any $i \in I_{G}$.

\rem{MPr} Let us note that a morphism between two normal crystals is automatically strict.
\erem

Let us say that $\mathbf{\iota}: \bB_{1} \ra \bB_{2}$ is an embedding of crystals if $\mathbf{\iota}$ is an injective morphism of crystals.

Let us say that $\bp: \bB_{1} \ra \bB_{2}$ is a retraction of crystals if $\bp$ is a morphism of crystals, $\bB_{2} \subset \on{Im}(\bp)$ and the restriction of $\bp$ to $\bp^{-1}(\bB_{2})$ is an isomorphism of crystals.

\rem{MPr} \label{1235} Let us note that if $\bp: \bB_{1} \ra \bB_{2}$ is a retraction of crystals then $\bB_{2} \simeq \bp^{-1}(\bB_{2}) \subset \bB_{1}.$ Thus every retraction of crystals $\bp: \bB_{1} \ra \bB_{2}$ induces the injection of crystals $\mathbf{\iota}: \bB_{2} \hookrightarrow \bB_{1}$.
\erem

\ssec{HWC}{Highest weight crystals} Let $\bB$ be a crystal. We say that $\bB$ is a highest weight crystal of weight $\la \in \La_{G}$ if there exists $\bb_{\la} \in \bB$, such that

(a) $wt(\bb_{\la}) = \la.$

(b) $\on{e}_{i}(\bb_{\la}) = 0$ for any $i \in I_{G}.$

(c) $\bB$ is generated by $\bb_{\la}.$

It is easy to see that if $\bB$ is a normal highest weight crystal of weight $\la$ then $\la \in \Lap_{G}$.

\ssec{CFC}{Closed families of $\check{G}$-crystals} Let us assume that for every $\la \in \Lap_{G}$ we are given a normal crystal $\bB^{G}(\la)$ of the highest weight $\la$. We say that the crystals $\bB^{G}(\la)$ form a closed family of $\check{G}$-crystals if for every $\la_{1}, \la_{2} \in \Lap_{G}$ there exists an embedding of crystals $\mathbf{\iota}_{\la_{1},\la_{2}}: \bB^{G}(\la_{1}+\la_{2}) \hookrightarrow \bB^{G}(\la_{1}) \otimes \bB^{G}(\la_{2})$.

\th{CFC}(see \cite[6.4.21]{Jo}) \label{1234} Let us assume that $G$ is of adjoint type. Then there exists the unique closed family of $\check{G}$-crystals.
\eth

\rem{CFC}
Let us point out that the condition for $G$ to be of adjoint type appears because of the following reason: $G$ is of adjoint type iff $\check{G}$ is simply connected. If $\check{G}$ is simply connected then it follows that any finite dimensional representation of $\check{\mathfrak{g}}$ integrates to the representation of $\check{G}$ so any closed family of crystals $\bB^{G}(\lambda)$ gives us the closed family of crystals for $\check{\mathfrak{g}}$ (as in \cite[Section 6.4.21]{Jo}).
\erem

Let us denote the family from Theorem \ref{1234} by $\bB^{\mathfrak{g}}(\la)$. Let us point out that from Remark \ref{1235} it follows that to construct $\bB^{\mathfrak{g}}(\la)$ it is enough to define normal crystals $\bB^{G}(\la)$ of the highest weight $\la$ for every $\lambda \in \Lap_{G}$ and retractions $\bp_{\la_{1},\la_{2}}: \bB^{G}(\la_{1}) \otimes \bB^{G}(\la_{2}) \ra \bB^{G}(\la_{1} + \la_{2}) \cup \{0\}$ for every $\la_{1}, \la_{2} \in \Lap_{G}$.

\quad

The main goal of this paper is to describe the family $\bB^{\mathfrak{g}}(\la)$ as sets of irreducible components of certain subvarieties in generalized transversal slices and to construct maps $\bp_{\la_{1},\la_{2}}$ using multiplication morphisms between various slices.

\sec{theo}{Main theorems} \label{42}

Let $\la \in \Lap_{G}, \mu \in \La_{G}$. Let $\bB^{G}(\la)_{\mu}$ denote the set of irreducible components of maximal dimension in $\ol{R}^{\la}_{\mu}$. Let $\bB^{G}(\la) := \coprod \limits_{\mu \in \La_{G}} \bB^{G}(\la)_{\mu}$.

\th{theo}
\label{7}
  (1) The map $\on{\bf{p}}^{\la}_{\mu}$ restricted to $\ol{R}^{\la}_{\mu}$ induces an isomorphism $\on{\bf{rp}}^{\la}_{\mu}: \ol{R}^{\la}_{\mu} \tilde{\ra} \ol{T}^{\la}_{\mu}$.

  (2) The set $\bB^{G}(\la)$ has a crystal structure of the highest weight $\la$. The collection $\{ \bB^{G}(\la) \}$ forms a closed family of crystals and the retraction morphisms $\bp_{\la_{1},\la_{2}}:\bB^{G}(\la_{1}) \otimes \bB^{G}(\la_{2}) \ra \bB^{G}(\la_{1}+\la_{2})\cup\{0\}$ are induced by the multiplication morphisms $\ol{\bf{\kappa}}^{\la_{1}, \la_{2}}_{\mu_{1}, \mu_{2}}:\ol{R}^{\la_{1}}_{\mu_{1}} \times \ol{R}^{\la_{2}}_{\mu_{2}} \ra \ol{R}^{\la_{1}+\la_{2}}_{\mu_{1}+\mu_{2}}$.
\eth

The proof of Theorem ~\ref{7} (1) will be given in Section ~\ref{22}, the proof of Theorem ~\ref{7}~(2) will be given in Section ~\ref{23}.

Let us give the examples of varieties $\ol{\CW}^{\la}_{\mu}, \ol{R}^{\la}_{\mu}$ in the case $G=GL_{2}$.

Let us fix $N,m \in \bN$ such that $m \leq N$. Let us suppose that $\lambda = N \cdot w_{1} = N \cdot \varepsilon_{1}$. Let us suppose that $\mu = (N-m) \cdot \varepsilon_{1} + m \cdot \varepsilon_{2}$.

It follows from \cite[Subsection 2(xii)]{BFN}] that the variety $\ol{\CW}^{\la}_{\mu}$ can be identified with the space of matrices $M:=\begin{bmatrix}
    D      & B  \\
    C      & A
\end{bmatrix} \in \on{Mat}_{2 \times 2}[z]$ such that A is a monic polynomial of degree m, degrees of $C,B$ are strictly less then $m$ and $\on{det}(M) = z^{N}$. Now it follows that

$$\ol{R}^{\la}_{\mu}=\{\begin{bmatrix}
    z^{N-m}     & 0  \\
    C      & z^{m}\end{bmatrix} | \on{deg}(C) \leq m-1\} \simeq \BA^{m}.$$

\rem{L1} Let us note that the condition $m \leq N$ precisely means that $\mu$ is a weight of $V^{\la}_{\check{G}}$ (compare with Lemma \ref{100001}).

\quad

Let us now suppose that $G=PGL_{3}$. Let us suppose that $\la = w_{1}, \mu = -w_{2}$ where $w_{1}, w_{2}$ are fundamental coweights of $G$. Let us note that $\mu = w_{0}(\la)$ thus it follows from \cite[See equality (3.6) in the proof of Theorem 3.2 ]{MV} that $\ol{T}^{\la}_{\mu}$ is an open cell inside $\ol{\on{Gr}}^{\la}_{G}$ so it is isomorphic to $\BA^{2}$. Now it follows from Theorem \ref{7} (1) that $\ol{R}^{\la}_{\mu} \simeq \BA^{2}$.
\erem

Let $P_{-} \subset G$ be a parabolic subgroup containing $B_{-}$.
Let $L\subset P_-$ be the corresponding Levi subgroup. Recall that $2{\rho}_{G,L}$ is the "subregular" dominant (i.e. the corresponding element of $\La_{G}$ lies inside $\Lap_{G}$ and $\langle 2{\rho}_{G,L}, \al_{i} \rangle > 0$ for any $i \notin I_{L}$) cocharacter of the torus $Z(L)$ (center of $L$). Let $\ol{R}^{\la}_{\mu,G,L}$ be the repellent of the action of $2{\rho}_{G,L}$ on $\ol{\CW}^{\la}_{\mu}$.

The slice $\ol{\CW}^{\la}_{\mu}$ is an affine scheme, thus the natural repelling morphism $\pi^{\la}_{\mu,G,L}:\ol{R}^{\la}_{\mu,G,L} \twoheadrightarrow (\ol{\CW}^{\la}_{\mu})^{2{\rho}_{G,L}}$ is a well defined algebraic morphism. Also we have the natural closed embedding $\iota^{\la}_{\mu,G,L}:\ol{R}^{\la}_{\mu,G,L} \hookrightarrow \ol{\CW}^{\la}_{\mu}$.

Let $\la-\mu=\sum\limits_{i \in I_{G}} n_{i}\alpha_{i}$. Let $\tilde{\la}:=\mu+\sum\limits_{i \in I_{L}} n_{i}\alpha_{i}$.

In Section ~\ref{24} we will construct a closed embedding $\varkappa^{\la}_{\mu,G,L}:(\ol{\CW}^{\la}_{\mu})^{2{\rho}_{G,L}} \hookrightarrow \ol{\CW}^{\tilde{\la}}_{\mu,L}$.

Let $\on{Res}^{\check{G}}_{\check{L}}: \on{Rep}(\check{G})\ra \on{Rep}(\check{L})$ denote the restriction functor from the category of finite dimensional $\check{G}$-modules to the category of $\check{L}$-modules. Our aim is to describe this functor in terms of generalized transversal slices. More concretely let $V_{\check{G}}^{\la}$ be an irreducible $\check{G}$-module of the highest weight $\la$. Let $(V^{\la}_{\check{G}})_{\mu}$ denote the $\check{T}$-weight $\mu$ component of $V_{\check{G}}^{\la}$. The representation $\on{Res}^{\check{G}}_{\check{L}}(V_{\check{G}}^{\la})$ can be decomposed into the sum of irreducible $\check{L}$-modules with some multiplicities. This decomposition is compatible with the $\check{T}$-weight decomposition. So we get the decomposition of $(V^{\la}_{\check{G}})_{\mu}$ into the direct sum of $\mu$-weight spaces of some irreducible $\check{L}$-modules $V^{\nu}_{\check{L}}$ with multiplicities.

Let $\La^{\la,+}_{\mu,G,L}$ denote the subset of $\nu \in \Lap_{L}$ such that the $\check{T}$-weight $\mu$ appears with nonzero multiplicity in the irreducible representation $V^{\nu}_{\check{L}}$ of $\check{L}$ and $V^{\nu}_{\check{L}}$ appears with nonzero multiplicity $n_{\nu}$ in the decomposition of $\on{Res}^{\check{G}}_{\check{L}}(V_{G}^{\la})$ into direct sum of irreducible $\check{L}$-modules.

$$(V^{\la}_{\check{G}})_{\mu}=\bigoplus\limits_{\nu \in \La^{\la,+}_{\mu,G,L}}{(V^{\nu}_{\check{L}})_{\mu}}^{\oplus{n_{\nu}}}.$$

For $\la \in \Lap_{G}$, $\mu \in \La_{G}$, let $\on{IC}(\ol{\CW}^{\la}_{\mu})$ be the $\on{IC}$ sheaf of $\ol{\CW}^{\la}_{\mu}$. For $\nu \in \Lap_{L}$ let $\on{IC}(\ol{\CW}^{\nu}_{\mu,L})$ be the $\on{IC}$ sheaf of $\ol{\CW}^{\nu}_{\mu,L}$.

 Consider the following diagram:
 \begin{equation}
\xymatrix{ \ol{\CW}^{\la}_{\mu} & \ar[l]_{\iota^{\la}_{\mu,G,L}} \ol{R}^{\la}_{\mu,G,L} \ar[r]^{\pi^{\la}_{\mu,G,L}} & (\ol{\CW}^{\la}_{\mu})^{2{\rho}_{G,L}}  \ar[r]^{\varkappa^{\la}_{\mu,G,L}} & \ol{\CW}^{\tilde{\la}}_{\mu,L}
}
\end{equation}
\th{theo}
\label{2}
The complex $(\varkappa^{\la}_{\mu,G,L} \circ \pi^{\la}_{\mu,G,L})_{*}(\iota^{\la}_{\mu,G,L})^{!}\on{IC}(\ol{\CW}^{\la}_{\mu})$ is perverse and is isomorphic to $\bigoplus\limits_{\nu \in \La^{+,\la}_{\mu,G,L}}(\ell^{\nu,\tilde{\la}}_{\mu,L})_{*}\on{IC}(\ol{\CW}^{\nu}_{\mu, L})^{\oplus{n_{\nu}}}.$
\eth

The proof will be given in Section ~\ref{24}.

\sec{L1}{Proof of Theorem ~\ref{7} (1)} \label{22}

\ssec{Mat}{Matrix description of slices} \label{21}
Recall an isomorphism (on the level of $\BC$-points, see \cite[2($\on{xi}$)]{BFN}) $$\Psi: \ol{\CW}^{\la}_{\mu} \simeq B[[z^{-1}]]_{1}z^{\mu}B_{-}[[z^{-1}]]_{1}\bigcap \ol{G[z]z^{\la}G[z]}.$$

\lem{L1} \label{25} The isomorphism $\Psi$ extends to the isomorphism between corresponding reduced ind-schemes (where $B[[z^{-1}]]_{1}z^{\mu}B_{-}[[z^{-1}]]_{1}\bigcap \ol{G[z]z^{\la}G[z]}$ is considered as a reduced locally closed ind-subscheme in $G(z)$).
\elem

\prf
First of all let us note that from \cite{BFN} Section 2(xi) it follows that the map $\Psi$ (considered as the map between $\ol{\CW}^{\la}_{\mu}$ and $G(z)$) is a locally closed embedding of ind-schemes.

Let $\bar{\iota}: \ol{B[[z^{-1}]]_{1}z^{\mu}B_{-}[[z^{-1}]]_{1}}\bigcap \ol{G[z]z^{\la}G[z]} \hookrightarrow G(z)$ denote the closed embedding of $\ol{B[[z^{-1}]]_{1}z^{\mu}B_{-}[[z^{-1}]]_{1}}\bigcap \ol{G[z]z^{\la}G[z]}$ to $G(z)$. Let us consider the following cartesian square:

\begin{equation}
\xymatrix{ P  \ar[r] \ar[d]^{\tilde{\bar{\iota}}} & \ol{B[[z^{-1}]]_{1}z^{\mu}B_{-}[[z^{-1}]]_{1}}\bigcap \ol{G[z]z^{\la}G[z]} \ar[d]^{\bar{\iota}} \\ \ol{\CW}^{\la}_{\mu}
 \ar[r] & G(z)
}
\end{equation}

where $P$ is the cartesian product $\ol{\CW}^{\la}_{\mu} \times_{G(z)} \ol{B[[z^{-1}]]_{1}z^{\mu}B_{-}[[z^{-1}]]_{1}}\bigcap \ol{G[z]z^{\la}G[z]}$.

Map $\tilde{\bar{\iota}}$ is the closed embedding which becomes an isomorphism on the level of $\BC$-points thus $\tilde{\bar{\iota}}$ is an isomorphism. So we get a morphism from $\ol{\CW}^{\la}_{\mu}$ to $\ol{B[[z^{-1}]]_{1}z^{\mu}B_{-}[[z^{-1}]]_{1}}\bigcap \ol{G[z]z^{\la}G[z]}$.

Now let us consider the following cartesian diagram (with the natural maps):

\begin{equation}
\xymatrix{ Q  \ar[r] \ar[d] &  \ar[d] B[[z^{-1}]]_{1}z^{\mu}B_{-}[[z^{-1}]]_{1} \bigcap \ol{G[z]z^{\la}G[z]} \\ \ol{\CW}^{\la}_{\mu}
 \ar[r] & \ol{B[[z^{-1}]]_{1}z^{\mu}B_{-}[[z^{-1}]]_{1}}\bigcap \ol{G[z]z^{\la}G[z]}
}
\end{equation}

where $Q:=\ol{\CW}^{\la}_{\mu} \times_{\ol{B[[z^{-1}]]_{1}z^{\mu}B_{-}[[z^{-1}]]_{1}}\bigcap \ol{G[z]z^{\la}G[z]}}  B[[z^{-1}]]_{1}z^{\mu}B_{-}[[z^{-1}]]_{1}\bigcap \ol{G[z]z^{\la}G[z]}$.

Note that the map from $B[[z^{-1}]]_{1}z^{\mu}B_{-}[[z^{-1}]]_{1}\bigcap \ol{G[z]z^{\la}G[z]}$ to $\ol{B[[z^{-1}]]_{1}z^{\mu}B_{-}[[z^{-1}]]_{1}}\bigcap \ol{G[z]z^{\la}G[z]}$ is the open embedding thus the map between $Q$ and $\ol{\CW}^{\la}_{\mu}$ is an open embedding which becomes an isomorphism on the level of $\BC$-points thus it is an isomorphism. So we have constructed the locally closed embedding $\Psi:\ol{\CW}^{\la}_{\mu} \ra B[[z^{-1}]]_{1}z^{\mu}B_{-}[[z^{-1}]]_{1}\bigcap \ol{G[z]z^{\la}G[z]}$.

One can see that the closure of $\ol{\CW}^{\la}_{\mu}$ inside $B[[z^{-1}]]_{1}z^{\mu}B_{-}[[z^{-1}]]_{1}\bigcap \ol{G[z]z^{\la}G[z]}$ coincides with $B[[z^{-1}]]_{1}z^{\mu}B_{-}[[z^{-1}]]_{1}\bigcap \ol{G[z]z^{\la}G[z]}$ (otherwise $\Psi$ could not be the isomorphism on the level of $\BC$-points). Thus the map $\Psi$ is the open embedding that is an isomorphism on the level of $\BC$-points so it is an isomorphism.
\epr

\rem{L1} The other way to prove that $\Psi$ is an isomorphism is to construct the inverse morphism $\Psi^{-1}$.

We construct $\Psi^{-1}$ as a morphism between the corresponding functors of points. Let $S$ be a test scheme. Take $\phi \in (B[[z^{-1}]]_{1}z^{\mu}B_{-}[[z^{-1}]]_{1}\bigcap \ol{G[z]z^{\la}G[z]})(S)$ i.e. a map from $S$ to $B[[z^{-1}]]_{1}z^{\mu}B_{-}[[z^{-1}]]_{1}\bigcap \ol{G[z]z^{\la}G[z]}$. It gives us a transition function of some $G$-bundle on $\BP^{1} \times S$. As in the end of \cite{BFN} Section 2(xi) it was observed we see that this bundle actually gives us a point in $\ol{\CW}^{\la}_{\mu}(S)$.

Thus we have constructed two algebraic morphisms between reduced ind-schemes $\ol{\CW}^{\la}_{\mu}$ and $B[[z^{-1}]]_{1}z^{\mu}B_{-}[[z^{-1}]]_{1}\bigcap \ol{G[z]z^{\la}G[z]}$ which are mutually inverse on the level of $\BC$-points. Thus they are isomorphisms (of reduced ind-schemes).
\erem

\rem{L1}
From Lemma ~\ref{25} it follows that $B[[z^{-1}]]_{1}z^{\mu}B_{-}[[z^{-1}]]_{1}\bigcap \ol{G[z]z^{\la}G[z]}$ is actually a scheme.
\erem

\lem{L1} Morphism $\Psi$ restricts to the isomorphism between $\ol{R}^{\la}_{\mu}$ and $z^{\mu}U_{-}[z^{-1}]_{1}\bigcap \ol{G[z]z^{\la}G[z]}$.
\elem

\prf
The Cartan torus $T$ acts on the scheme $B[[z^{-1}]]_{1}z^{\mu}B_{-}[[z^{-1}]]_{1}\bigcap \ol{G[z]z^{\la}G[z]}$ via conjugation. So the repellent to the point $z^{\mu}$ is $T[[z^{-1}]]_{1}z^{\mu}U_{-}[[z^{-1}]]_{1} \bigcap \ol{G[z]z^{\la}G[z]}$. Note that $\ol{G[z]z^{\la}G[z]}\subset G(z)$ so $T[[z^{-1}]]_{1}z^{\mu}U_{-}[[z^{-1}]]_{1} \bigcap  \ol{G[z]z^{\la}G[z]} = T[z^{-1}]_{1}z^{\mu}U_{-}[z^{-1}]_{1} \bigcap  \ol{G[z]z^{\la}G[z]}$. But the scheme $T[z^{-1}]_{1}=1$, so $\ol{R}^{\la}_{\mu} = z^{\mu}U_{-}[z^{-1}]_{1}\bigcap \ol{G[z]z^{\la}G[z]}$.
\epr

\lem{L1} The natural map $\xi : U_{-}[z^{-1}]_{1} \ra U_{-}(\CK)/U_{-}(\CO)$ (induced by the embedding $U_{-}[z^{-1}]_{1} \hookrightarrow U_{-}(\CK)$) is an isomorphism of ind-schemes.
\elem

\prf
Let us show that the natural morphism $\on{\bf{mult}}: U_{-}[z^{-1}]_{1} \times U_{-}[[z]] \ra U_{-}((z))$ is an isomorphism of ind-schemes. In other words we want to show that for any test $\BC$-algebra $R$ the morphism $\on{\bf{mult}}: U_{-}(R[z^{-1}]_{1}) \times U_{-}(R[[z]]) \ra U_{-}(R((z)))$ is an isomorphism.
The unipotent group-scheme $U_{-}$ can be filtered by normal subgroups with a successive quotients isomorphic to $\BG_{a}$. Thus it is enough to prove our statement for $\BG_{a}$. For the group-scheme $\BG_{a}$ we have $\BG_{a}(R[z^{-1}])_{1}=z^{-1}R[z^{-1}]$, $\BG_{a}(R[[z]])=R[[z]]$, $\BG_{a}(R((z)))=R((z))$. Obviously the addition map $z^{-1}R[z^{-1}] \times R[[z]] \ra R((z))$ is an isomorphism.

We obtain the isomorphism ${\bf{mult}}$ between two right $U_{-}[[z]]$-torsors. It induces the desired isomorphism $\xi$.
\epr

\lem{L1} \label{60} The morphism $\on{\bf{rp}}^{\la}_{\mu}$ is a locally closed embedding.
\elem

\prf Let us consider the locally closed embedding $\on{\bp}^{\la}_{\mu}\times\on{\bq}^{\la}_{\mu}:\ol{\CW}^{\la}_{\mu} \hookrightarrow \ol{\on{Gr}}^{\la}_G \times Z^{-w_{0}(\la-\mu)}$. It gives us the locally closed embedding $\on{\bf{rp}}^{\la}_{\mu}\times\on{\bf{rq}}^{\la}_{\mu}:\ol{R}^{\la}_{\mu} \hookrightarrow \ol{T}^{\la}_{\mu} \times RZ^{-w_{0}(\la-\mu)}$ where ${RZ}^{-w_{0}(\la-\mu)}$ is the repellent of $Z^{-w_{0}(\la-\mu)}$ to the fixed point $\bq^{\la}_{\mu}(z^{\mu})$. But ${RZ}^{-w_{0}(\la-\mu)}$ coincides with $\bq^{\la}_{\mu}(z^{\mu})$ (Actually, $\on{\bf{rq}}^{\la}_{\mu}$ maps $\ol{R}^{\la}_{\mu}$ to the central fiber of $Z^{-w_{0}(\la-\mu)}$ under the factorization morphism and it is known \cite[Theorem 2.7]{BFK} (also see~\cite[Section 1]{M} for the analogous statement for the whole Zastava space) to be isomorphic to $\ol{T}_{w_{0}(\la-\mu)}\cap S_{0}$; in particular, it lies inside the {\em{attractor}} $S_{0}$. See also Corollary ~\ref{45} below for the proof of more general statement). Thus $\on{\bf{rp}}^{\la}_{\mu}\times\on{\bf{rq}}^{\la}_{\mu}$ induces the locally closed embedding from $\ol{R}^{\la}_{\mu}$ to $\ol{T}^{\la}_{\mu} \times \on{\bq}^{\la}_{\mu}(z^{\mu})$ so $\on{\bf{rp}}^{\la}_{\mu}$ is the locally closed embedding.
\epr

\cor{L1} From Lemma ~\ref{60} it follows that to prove Theorem ~\ref{7} (1) it is enough to show that $\on{\bf{rp}}^{\la}_{\mu}$ is surjective.
\ecor

\ssec{Pr}{First proof of Theorem ~\ref{7} (1)}
\prf Let us construct the inverse map $({\on{\bf{rp}}^{\la}_{\mu}})^{-1} : (U_{-}(\CK)z^{\mu}/U_{-}(\CO)) \cap \ol{\on{Gr}}^{\la}_G \ra z^{\mu}U_{-}[z^{-1}]_{1}\bigcap \ol{G[z]z^{\la}G[z]}$. Let us take a point $[u_{-}z^{\mu}] \in \ol{T}^{\la}_{\mu}$ and rewrite it as $[z^{\mu}\tilde{u}_{-}]$. Let $({\on{\bf{rp}}^{\la}_{\mu}})^{-1}(u_{-}z^{\mu}) := z^{\mu}\xi^{-1}(\tilde{u}_{-})$.  We have to check that $z^{\mu}\xi^{-1}(\tilde{u}_{-}) \in z^{\mu}U_{-}[z^{-1}]_{1}\bigcap \ol{G[z]z^{\la}G[z]}$. The point $\xi^{-1}(\tilde{u}_{-})$ lies in $U_{-}[z^{-1}]_{1},$ so  $z^{\mu}\xi^{-1}(\tilde{u}_{-}) \in z^{\mu}U_{-}[z^{-1}]_{1}$. The point $z^{\mu}\xi^{-1}(\tilde{u}_{-})$ considered as a point of $\on{Gr}_G$ coincides with $u_{-}z^{\mu}$ so the corresponding trivialization has degree $\leq \la$, thus $z^{\mu}\xi^{-1}(\tilde{u}_{-}) \in \ol{G[z]z^{\la}G[z]}$. It is easy to see that the maps ${\on{\bf{rp}}^{\la}_{\mu}}$ and $({\on{\bf{rp}}^{\la}_{\mu}})^{-1}$ are mutually inverse.
\epr

\ssec{Pr}{Second proof of Theorem ~\ref{7} (1)} \label{20}
\prf

The first proof of Theorem ~\ref{7} (1) used the "matrix" (see Lemma ~\ref{25}) description of slices. The advantage of the second proof is that it does not.

We have already proved that the map $\on{\bf{rp}}^{\la}_{\mu}$ is a locally closed embedding. So we have to show that it is surjective on the level of $\BC$-points. Let us construct a section $s^{\la}_{\mu}: \ol{T}^{\la}_{\mu}(\BC) \ra  \ol{R}^{\la}_{G,L, \mu}(\BC)$ of the map $\on{\bf{rp}}^{\la}_{\mu}$.

By the definition we have the following cartesian square:

\begin{equation}
\xymatrix{ \ol{T}^{\la}_{\mu}(\BC)  \ar[r] \ar[d] &  \ol{\on{Gr}}^{\la}_{G}(S) \ar[d] \\
\on{Gr}_{B_{-},\mu}(\BC) \ar[r]& \on{Gr}_{G}(\BC)
}
\end{equation}

Note that $\on{Gr}_{B_{-},\mu}$ is isomorphic to $\on{Gr}_{w_{0}(\mu),B}$ via conjugation by $\tilde{w_{0}} \in G$ (some fixed representative of $w_{0}$) inside $\on{Gr}_{G}$.
So $\ol{T}^{\la}_{\mu}(\BC)$ is the following moduli space:

(a) A $G$-bundle $\CP$ on $\BP^{1}$.

(b) A trivialization $\sig: \CP_{triv}|_{(\BP^{1} \sm 0)} \tilde{\ra} \CP|_{(\BP^{1} \sm 0)}$ such that the corresponding point of $\on{Gr}_{G}(\BC)$ lies in $\ol{\on{Gr}}^{\la}_{G}(\BC)$.

(c) A $B$-bundle $\CP^{B}$ of degree $w_{0}(\mu)$ with a trivialization $'\sigma^{B}: \CP^{B}_{triv}|_{(\BP^{1} \sm 0)} \tilde{\ra} \CP^{B}|_{(\BP^{1} \sm 0)}$.

(d) An isomorphism $\psi$ between $G\times_{B} \CP^{B}$ and $\CP$ such that the corresponding trivializations differ by the action of $\tilde{w_{0}}$ (that means in particular that in the fibers of $\CP$ over $\infty $ our $B$-structure is identified with $B_{-}$ via $\sig$).

Thus we get the point $(\CP,\sig,\CP^{B})$ of $\ol{\CW}^{\la}_{\mu}(\BC)$. It is easy to see that it lies inside $\ol{R}^{\la}_{\mu}(\BC)$.

Thus we have constructed a map $s^{\la}_{\mu}:  \ol{T}^{\la}_{\mu}(\BC) \ra  \ol{R}^{\la}_{G,L, \mu}(\BC)$.

Directly from the definitions $\on{\bf{rp}}^{\la}_{\mu} \circ s^{\la}_{\mu}=\on{Id}$.
\epr

\rem{pr} Note that the isomorphism between $\on{Gr}_{B_{-},\mu}$ and $\on{Gr}_{B,w_{0}(\mu)}$ is not canonical (it depends on the choice of a representative of $w_{0}$ inside $G$) while it turns out that the morphism $s^{\la}_{\mu}$ does not depend on this choice.
\erem

\sec{Res}{Proof of Theorem ~\ref{2}} \label{24}

Let $\la-\mu=\sum\limits_{i \in I_{G}} n_{i}\alpha_{i}$, $\tilde{\la}:=\mu+\sum\limits_{i \in I_{L}} n_{i}\alpha_{i}$.
Let $\La^{\la,+}_{G,L}$ denote the subset of $\La^{+}_{L}$ consisting of $\nu$ such that $V^{\nu}_{\check{L}}$ appears with nonzero multiplicity in the decomposition of $\on{Res}^{\check{G}}_{\check{L}}(V_{\check{G}}^{\la})$ into direct sum of irreducible $\check{L}$-modules.

\lem{Res} \label{15}
$(\ol{\on{Gr}}^{\la}_{G})^{2{\rho}_{G,L}}=\bigcup\limits_{\nu \in \La^{\la,+}_{G,L}}\ol{\on{Gr}}^{\nu}_{L}$.
\elem

\prf
First of all let us prove that $\bigcup\limits_{\nu \in \La^{\la,+}_{G,L}}\ol{\on{Gr}}^{\nu}_{L} \subset (\ol{\on{Gr}}^{\la}_{G})^{2{\rho}_{G,L}}$. The subvariety $(\ol{\on{Gr}}^{\la}_{G})^{2{\rho}_{G,L}}$ is closed inside $\on{Gr}_{G}$ and $L(\CO)$-invariant. So it is enough to show that for any $\nu \in \La^{\la,+}_{G,L}$, the point
$z^{\nu}$ lies in $(\ol{\on{Gr}}^{\la}_{G})^{2{\rho}_{G,L}}$. The representation $V^{\nu}_{\check{L}}$ appears with nonzero multiplicity in the decomposition of $\on{Res}^{\check{G}}_{\check{L}}(V_{G}^{\la})$ into the direct sum of irreducible $\check{L}$-modules, thus the weight $\nu$ appears with nonzero multiplicity in the decomposition of $\on{Res}^{\check{G}}_{\check{T}}(V_{G}^{\la})$ into the direct sum of irreducible $\check{T}$-modules, so that $z^{\nu} \in \ol{\on{Gr}}^{\la}_{G}$. Thus $z^{\nu} \in (\ol{\on{Gr}}^{\la}_{G})^{2{\rho}_{G,L}}$.

Now let us prove that $(\ol{\on{Gr}}^{\la}_{G})^{2{\rho}_{G,L}} \subset \bigcup\limits_{\nu \in \La^{\la,+}_{G,L}}\ol{\on{Gr}}^{\nu}_{L}$. It is enough to show that if $\on{Gr}^{'\nu}_{L}\cap(\ol{\on{Gr}}^{\la}_{G})^{2{\rho}_{G,L}}$ is nonempty for some $'\nu \in \Lap_{L}$ then there exists $\nu \in \La^{\la,+}_{G,L}$ such that $'\nu$ appears in $V^{\nu}_{\check{L}}$ with a nonzero multiplicity.

Let us suppose that $\on{Gr}^{'\nu}_{L}\cap(\ol{\on{Gr}}^{\la}_{G})^{2{\rho}_{G,L}}$ is nonempty. Let $'\nu^{\on{dom}}$ be the $G$-dominant cocharacter in the $W_{G}$-orbit of $'\nu$. One can note that $z^{'\nu} \in (\on{Gr}^{'\nu^{\on{dom}}}_{G})^{2{\rho}_{G,L}}$ thus $\on{Gr}^{'\nu}_{L} \subset (\on{Gr}^{'\nu^{\on{dom}}}_{G})^{2{\rho}_{G,L}}$. So $(\on{Gr}^{'\nu^{\on{dom}}}_{G}) \cap \ol{\on{Gr}}^{\la}_{G}$ is nonempty, thus $'\nu^{\on{dom}} \leq \la$. So $'\nu$ appears as a $\check{T}$-weight of $V^{\la}_{\check{G}}$, thus $'\nu$ appears as a $\check{T}$-weight of $V^{\nu}_{\check{L}}$ for some $\nu \in \La^{\la,+}_{G,L}.$
\epr

\rem{Res} Lemma ~\ref{15} is a generalization of the following statement: for $\mu \in \La_{G}$ $\la \in \Lap_{G}$, $z^{\mu} \in \ol{\on{Gr}}^{\la}_{G}$ iff $\mu$ is a $\check{T}$-weight of the representation $V^{\la}_{\check{G}}$.
\erem

Let $\theta \in \La_{G,P_{-}}$. Let $\La^{\la,+}_{\theta,G,L}$ denote the subset of $\Lap_{L}$ consisting of all $\nu$ such that $V^{\nu}_{\check{L}}$ appears with nonzero multiplicity in the decomposition of $\on{Res}^{\check{G}}_{\check{L}}(V^{\la}_{\check{G}})$
and $Z(\check{L})$ acts on $V^{\nu}_{\check{L}}$ via the character $\theta$. For $\mu \in \La_{G}$ let $\bar{\mu}:= \al_{G,P_{-}}(\theta)$.

\lem{Res}
\label{3}
$(\ol{\on{Gr}}^{\la}_{G} \cap \on{Gr}_{P_{-},\bar{\mu}})^{2{\rho}_{G,L}}=\bigcup\limits_{\nu \in \La^{\la,+}_{\bar{\mu},G,L}}\ol{\on{Gr}}^{\nu}_{L} \subset \ol{\on{Gr}}^{\tilde{\la}}_{L}$.
\elem

\prf
$(\ol{\on{Gr}}^{\la}_{G} \cap \on{Gr}_{P_{-},\bar{\mu}})^{2{\rho}_{G,L}}=\bigcup \limits_{\nu \in \La^{\la,+}_{G,L}}\ol{\on{Gr}}^{\nu}_{L} \cap \on{Gr}_{L,\bar{\mu}}=\bigcup \limits_{\nu \in \La^{\la,+}_{\bar{\mu},G,L}}\ol{\on{Gr}}^{\nu}_{L}.$

The latter equality holds because of the following inclusion:\\
$\on{Gr}^{\nu}_{L} \subset \on{Gr}_{\bar{\nu}, L}$, where $\bar{\nu}:=\alpha_{G,P_{-}}(\nu)$. To see that let us note that $\on{Gr}^{\nu}_{L}$ is connected thus it is enough to show that the intersection $\on{Gr}^{\nu}_{L} \cap \on{Gr}_{\bar{\nu}, L}$ is nonempty. Obviously $z^{\nu} \in \on{Gr}^{\nu}_{L} \cap \on{Gr}_{\bar{\nu}, L}.$
\epr

\lem{Res} \label{17}
$(\ol{G[z]z^{\la}G[z]})^{2{\rho}_{G,L}}=\bigcup\limits_{\nu \in \La^{\la,+}_{G,L}}\ol{L[z]z^{\nu}L[z]}$.
\elem

\prf
Same argument as in Lemma ~\ref{15}.
\epr

Let us denote by $\bp^{\tilde{\la}}_{L,\mu}:\ol{\CW}^{\tilde{\la}}_{\mu,L}\ra \ol{\on{Gr}}^{\tilde{\la}}_{L}$ the natural "forgetting" morphism.
\lem{Res}
\label{4} The variety $({\ol{\CW}^{\la}_{\mu}})^{2{\rho}_{G,L}}$ is isomorphic to $(\on{\bp}^{\tilde{\la}}_{\mu,L})^{-1}(\bigcup\limits_{\nu \in \La^{\la,+}_{G,L}, \mu \leq_{L} \nu} \ol{\on{Gr}}^{\nu}_{L})=\bigcup\limits_{\nu \in \La^{\la,+}_{G,L}, \mu \leq_{L} \nu} \ol{\CW}^{\nu}_{\mu,L}$. In particular we have the following embeddings:

$$\bigcup\limits_{\nu \in \La^{\la,+}_{\mu,G,L}} \ol{\CW}^{\nu}_{\mu,L} \subset ({\ol{\CW}^{\la}_{\mu}})^{2{\rho}_{G,L}} \subset {\ol{\CW}^{\tilde{\la}}_{\mu,L}}.$$
\elem

\prf
$(\ol{\CW}^{\la}_{\mu})^{2{\rho}_{G,L}} \simeq$ $(B[[z^{-1}]]_{1}z^{\mu}B_{-}[[z^{-1}]]_{1})^{2{\rho}_{G,L}}\bigcap \ol{G[z]z^{\la}G[z]}^{2{\rho}_{G,L}} =$

$B_{L}[[z^{-1}]]_{1}z^{\mu}B_{L,-}[[z^{-1}]]_{1}\bigcap \ol{G[z]z^{\la}G[z]}^{2{\rho}_{G,L}}=$

$B_{L}[[z^{-1}]]_{1}z^{\mu}B_{L,-}[[z^{-1}]]_{1}\bigcap \bigcup\limits_{\nu \in \La^{\la,+}_{G,L}}\ol{L[z]z^{\nu}L[z]}=$
$\bigcup\limits_{\nu \in \La^{\la,+}_{G,L}}\ol{\CW}^{\nu}_{\mu,L}$

(second equality follows from Lemma ~\ref{17}).

Note that a variety $\ol{\CW}^{\nu}_{\mu,L}$ is nonempty iff $\mu \leq_{L} \nu$ thus the desired follows.
\epr

\rem{Res}
We will see below that both embeddings of Lemma~\ref{4} are not necessarily isomorphisms.
\erem

For $\la \in \Lap_{G}$, $\mu \in \La_{G}$, $\nu \in \La^{\la,+}_{\mu,G,L}$, let $\imath^{\nu}_{\mu,G,L}: \ol{\CW}^{\nu}_{\mu,L} \hookrightarrow (\ol{\CW}^{\la}_{\mu})^{2{\rho}_{G,L}}$ be the closed embedding of Lemma ~\ref{4}.

From Lemma ~\ref{4} it follows that we have a natural closed embedding $\varkappa^{\la}_{\mu,G,L}: (\ol{\CW}^{\la}_{\mu})^{2\check{\rho}_{G,L}} \hookrightarrow \ol{\CW}^{\tilde{\la}}_{\mu,L}$. This embedding is not always an isomorphism (but in the case of dominant $\mu$ it is obviously an isomorphism).

\ssec{first example}{Example} Let $G=\check{G}=GL_{4}(\BC)$, $L=GL_{2}(\BC)\times GL_{2}(\BC)$, $\varepsilon_{1}, \varepsilon_{2}, \varepsilon_{3}, \varepsilon_{4}$ are the natural generators of the character lattice $\La_{G}$. Let $\la:=2\varepsilon_{1}-2\varepsilon_{4}$, $\mu:=-\varepsilon_{2}+\varepsilon_{3}$.

So $\la-\mu=2(\varepsilon_{1}-\varepsilon_{2})+3(\varepsilon_{2}-\varepsilon_{3})+2(\varepsilon_{3}-\varepsilon_{4})$. Thus $\tilde{\la}=2(\varepsilon_{1}-\varepsilon_{2})+2(\varepsilon_{3}-\varepsilon_{4})+\mu=2\varepsilon_{1}-3\varepsilon_{2}+3\varepsilon_{3}-2\varepsilon_{4}$. We want to show that $\tilde{\la}$ is not a weight of irreducible representation of $GL_{4}(\BC)$ with the highest weight $\la$.

To check it we should consider the dominant character $\tilde{\la}^{\on{dom}}$ in the orbit of $\tilde{\la}$ under the action of the Weyl group of $GL_{4}(\BC)$ and to show that it is not less than or equal to $\la$. It is easy to see that $\tilde{\la}^{\on{dom}}=\la+(\varepsilon_{2}-\varepsilon_{3})>\la$.
Thus $\tilde{\la}^{\on{dom}}$ is not a weight of irreducible representation of $GL_{4}(\BC)$ with the highest weight $\la$. So $\tilde{\la}$ is not a weight of irreducible representation of $GL_{4}(\BC)$ with the highest weight $\la$.

\bigskip

From Lemma ~\ref{4} it follows that $({\ol{\CW}^{\la}_{\mu}})^{2\check{\rho}_{G,L}}=\bigcup\limits_{\nu \in \La^{\la,+}_{G,L}, \mu \leq_{L} \nu} \ol{\CW}^{\nu}_{\mu,L}$, thus we get a closed embedding $\bigcup\limits_{\nu \in \La^{\la,+}_{\mu,G,L}} \ol{\CW}^{\nu}_{\mu,L} \hookrightarrow ({\ol{\CW}^{\la}_{\mu}})^{2\check{\rho}_{G,L}}$. Let us give an example when this embedding is not an isomorphism.

\ssec{second example}{Example} We have to find $\la \in \Lap_{G}$, $\mu \in \La_{G}$ and $\nu \in \Lap_{L}$ such that $\mu \leq_{L} \nu$ and $\mu$ appears as a $\check{T}$-weight of $V^{\la}_{\check{G}}$ but $\mu$ does not appear as a $\check{T}$-weight of $V^{\nu}_{\check{L}}$.

Let $G=\check{G}=GL_{3}$, $L=GL_{2}\times GL_{1}$, $\la = 2\varepsilon_{1}+\varepsilon_{2}$, $\mu = 2\varepsilon_{2}+\varepsilon_{3}$, $\nu = \varepsilon_{1}+\varepsilon_{2}+\varepsilon_{3}$. Let us consider an irreducible representation $V^{\la}_{G}$ of $G$ with highest weight $\la$. Let us show that as an $L$-module it contains the one dimensional subrepresentation of weight $\nu$. To see that let us note that $V^{\la}_{G}$ appears as a "highest" subrepresentation of $G$-module $\BC^{3} \otimes \La^{2}(\BC^{3})$. Let $v_{1}, v_{2}, v_{3}$ be the natural basis in $\BC^{3}$. It is easy to see that one-dimensional vector space generated by $2v_{3} \otimes (v_{1} \wedge v_{2}) + v_{2} \otimes (v_{1} \wedge v_{3}) - v_{1} \otimes (v_{2} \wedge v_{3})$ lies in $V^{\la}_{G}$, is stable under the action of $L$ and has the weight $\nu$.

Note that $\mu <_{L} \nu$ and $\mu$ appears as a weight of $V_{G}^{\la}$ but does not appear as a $\check{T}$-weight of $V^{\nu}_{\check{L}}$.

Let us take $\al \in \La^{\on{pos}}_{G}$. Let us denote by $SZ^{\al}$, $RZ^{\al}$ attractor and repellent of the $\BC^{*}$-action on $Z^{\al}$ via $2{\rho}_{G,L}$.

\lem{Res}(see ~\cite[Theorem 2.7]{BFK}) \label{55}
$SZ^{\al}$ coincides with the whole space $Z^{\al}$.
\elem

\prf
Let us take $f \in Z^{\al}$. We have to show that there exists a limit $\on{lim}_{t \ra 0}2{\rho}_{G,L}(t)f \in Z^{\al}$. Let us consider the locally closed embedding $Z^{\al} \hookrightarrow \on{QMaps}(\BP^{1},\CB)$ of Zastava space to the space of quasi-maps of degree $\al$ from $\BP^{1}$ to the flag variety $\CB$. Note that $\on{QMaps}^{\al}(\BP^{1},\CB)$ is projective; thus the limit $\on{lim}_{t \ra 0}2{\rho}_{G,L}(t)f$ exists inside $\on{QMaps}^{\al}(\BP^{1},\CB)$. Let us denote it $f_{0}$. We have to show that $f_{0}$ actually lies inside $Z^{\al}$. To see
this we have to check that $f_{0}(\infty)=B_{-}$ and $f_{0}$ has no defect in $\infty$. The first property is obvious because $B_{-} \in \CB$ is fixed under the action of $2{\rho}_{G,L}$ and $f(\infty)=B_{-}$. To check the second property let us note that $f$ has no defect at $\infty$ and $f(\infty)=B_{-}$ thus there exists an open subset $V \subset \BP^{1}$ such that $f$ restricted to $V$ defines an actual map $f|_{V}: V \ra U \cdot B_{-}$ where $U \cdot B_{-}$ is an open Bruhat cell inside $\CB$.

Now the desired claim follows from the fact that for any $\tilde{f}: V \ra U \cdot B_{-}$ there exists the limit $\on{lim}_{t \ra 0}2{\rho}_{G,L}(t)\tilde{f}$ inside $\on{Maps}(V,U\cdot B_{-})$ (because $U\cdot B_{-}$ is contained in the attractor of $G$ under the action of $2{\rho}_{G,L}$).

\epr

\cor{Res} \label{45}
The natural closed embedding $(Z^{\al})^{2{\rho}_{G,L}} \hookrightarrow RZ^{\al}$ is
an isomorphism.
\ecor

\prf
From Lemma \label{55} it follows that $SZ^{\al}=Z^{\al}$, thus $RZ^{\al}=RZ^{\al} \cap SZ^{\al}$. Note that $Z^{\al}$ is an affine variety, so $RZ^{\al} \cap SZ^{\al} = (Z^{\al})^{2{\rho}_{G,L}}$.

\epr

\lem{Res} \label{6}
The following square is cartesian:

\begin{equation}
\xymatrix{ \ol{R}^{\la}_{G,L, \mu}  \ar[r] \ar[d]^{\pi^{\la}_{G,L,\mu}} &  \ol{\on{Gr}}^{\la}_{G} \cap \on{Gr}_{P_{-}, \bar{\mu}} \ar[d]^{\tilde{\pi}}\\
({\ol{\CW}^{\la}_{\mu}})^{2{\rho}_{G,L}} \ar[r]^{(\on{\bp}^{\la}_{\mu})^{2{\rho}_{G,L}}} & \ol{\on{Gr}}^{\la}_{G} \cap \on{Gr}_{L,\bar{\mu}}
}
\end{equation}

where $\bar{\mu}:=\alpha_{G,P_{-}}(\mu).$
\elem

\rem{Res}
 Lemma ~\ref{6} is a "relative" version of Theorem ~\ref{7} (1). We will give two proofs of Lemma ~\ref{6} below. They are respective generalizations of our two proofs of Theorem ~\ref{7} (1).
\erem

\ssec{Res}{First proof of Lemma ~\ref{6}}
\prf
We have to show that $\ol{R}^{\la}_{G,L, \mu}\simeq (\ol{\on{Gr}}^{\la}_{G} \cap \on{Gr}_{P_{-},\bar{\mu}})\times_{\ol{\on{Gr}}^{\la}_{G} \cap \on{Gr}_{L,\bar{\mu}}} ({\ol{\CW}^{\la}_{\mu}})^{2{\rho}_{G,L}}$. A morphism ${\bf{rp}}^{\la}_{G,L,\mu}$ from $\ol{R}^{\la}_{G,L, \mu}$ to $(\ol{\on{Gr}}^{\la}_{G} \cap \on{Gr}_{P_{-},\bar{\mu}})\times_{\ol{\on{Gr}^{\la}_{G}} \cap \on{Gr}_{L,\bar{\mu}}} ({\ol{\CW}^{\la}_{\mu}})^{2{\rho}_{G,L}}$ can be constructed in an obvious way.

It is easy to see that the morphism $\on{{\bf{rp}}}^{\la}_{G,L, \mu}$ is a locally closed embedding. It follows from Corollary ~\ref{45} using the same observations as in ~\ref{20}. Thus we just have to show that $\on{{\bf{rp}}}^{\la}_{G,L, \mu}$ is surjective.

Let us construct the inverse morphism. Let us take points $[p] \in \ol{\on{Gr}}^{\la}_{G} \cap \on{Gr}_{P_{-},\bar{\mu}}$, $l \in ({\ol{\CW}^{\la}_{\mu}})^{2{\rho}_{G,L}}$ such that $\tilde{\pi}([p])=\on{(\bp^{\la}_{\mu})^{2{\rho}_{G,L}}}(l)$. Using these data we want to construct a point in $\ol{R}^{\la}_{G,L,\mu}$. Clearly, $\tilde{\pi}([p])=\on{(\bp^{\la}_{\mu})^{2{\rho}_{G,L}}}(l)$ inside $\on{Gr}_{L,\bar{\mu}}$, so there exists a unique $\tilde{l} \in L(\CO)$ such that $p = l \tilde{l} \on{mod} U_{P_{-}}(\CK)$. Thus changing the representative of $[p] \in \on{Gr}_{P_{-},\bar{\mu}}$ from $p$ to $p\tilde{l}^{-1}$ we can assume that $p \in P_{-}((\CK))$ contracts via the action of $2{\rho}_{G,L}$ to the point $l$. Thus $p=lu_{-}$ for some $u_{-}\in U_{P_{-}}(\CK)$. Using the isomorphism $\xi_{U_{P_{-}}}$ between $U_{P_{-}}[z^{-1}]_{1}$ and $U_{P_{-}}(\CK)/U_{P_{-}}(\CO)$ we get the point $\xi_{U_{P_{-}}}^{-1}(u_{-}) \in U_{P_{-}}[z^{-1}]_{1}$. L $({\bf{rp}}^{\la}_{G,L,\mu})^{-1}([p],l):=l\xi_{U_{-}}^{-1}(u_{-})$. The point $({\bf{rp}}^{\la}_{G,L,\mu})^{-1}([p],l)$ lies in $\ol{R}^{\la}_{G,L, \mu}$ and $[(\bp^{\la}_{G,L,\mu})^{-1}([p],l)]=[p]$. To see that it lies in $\ol{R}^{\la}_{G,L, \mu}$ we have to show that $l\xi_{U_{P_{-}}}^{-1}(u_{-}) \in B_{L}[[z^{-1}]]_{1}z^{\mu}B_{-}[[z^{-1}]]_{1}\bigcap \ol{G[z]z^{\la}G[z]}$. The point $l$ lies in $B_{L}[[z^{-1}]]_{1}z^{\mu}B_{L, -}[[z^{-1}]]_{1}$ and $\xi_{U_{P_{-}}}^{-1}(u_{-}) \in U_{P_{-}}[z^{-1}]_{1}$, so $l\xi_{P_{-}}^{-1}(u_{-}) \in B_{L}[[z^{-1}]]_{1}z^{\mu}B_{-}[[z^{-1}]]_{1}$. To see that $l\xi_{P_{-}}^{-1}(u_{-}) \in \ol{G[z]z^{\la}G[z]}$, note that $[l\xi_{P_{-}}^{-1}(u_{-})]=[p]$ and $[p] \in \ol{\on{Gr}}^{\la}_{G}$. Thus the trivialization of the corresponding bundle has a pole of degree $\leq \la$, so $l\xi_{P_{-}}^{-1}(u_{-}) \in \ol{G[z]z^{\la}G[z]}$.
\epr

\ssec{Res}{Second proof of Lemma ~\ref{6} (sketch)}

We have to show that the following square is cartesian:
\begin{equation}
\xymatrix{ \ol{R}^{\la}_{G,L, \mu}  \ar[r] \ar[d]^{\pi^{\la}_{G,L,\mu}} &  \ol{\on{Gr}}^{\la}_{G} \cap \on{Gr}_{P_{-},\bar{\mu}} \ar[d]^{\tilde{\pi}}\\
({\ol{\CW}^{\la}_{\mu}})^{2{\rho}_{G,L}} \ar[r]^{(\on{\bp}^{\la}_{\mu})^{2{\rho}_{G,L}}} & \ol{\on{Gr}}^{\la}_{G} \cap \on{Gr}_{L,\bar{\mu}}
}
\end{equation}
\prf
Let us consider the natural morphism $\on{{\bf{rp}}}^{\la}_{G,L, \mu}:\ol{R}^{\la}_{G,L,\mu} \ra (\ol{\on{Gr}}^{\la}_{G} \cap \on{Gr}_{P_{-},\bar{\mu}})\times_{\ol{\on{Gr}}^{\la}_{G} \cap \on{Gr}_{L,\bar{\mu}}} ({\ol{\CW}^{\la}_{\mu}})^{2{\rho}_{G,L}}$. We already know that the morphism $\on{{\bf{rp}}}^{\la}_{G,L, \mu}$ is a locally closed embedding.

It is enough to construct a section $s^{\la}_{G,L,\mu}:(\ol{\on{Gr}}^{\la}_{G} \cap \on{Gr}_{P_{-},\bar{\mu}})\times_{\ol{\on{Gr}}^{\la}_{G} \cap \on{Gr}_{L,\bar{\mu}}} ({\ol{\CW}^{\la}_{\mu}})^{2{\rho}_{G,L}} \ra \ol{R}^{\la}_{G,L, \mu}$ of the map $\on{{\bf{rp}}}^{\la}_{G,L, \mu}$
on the level of $\BC$-points.
Note that $\on{Gr}_{P_{-},\bar{\mu}}$ is isomorphic to $\on{Gr}_{P,w_{0}({w_{0}}^{L})^{-1}(\bar{\mu})}$ inside $\on{Gr}_{G}$ where $w_{0}({w_{0}}^{L})^{-1}(\bar{\mu}):=\al_{G,P_{-}}(w_{0}({w_{0}}^{L})^{-1}(\mu))$.

A $\BC$-point of $(\ol{\on{Gr}}^{\la}_{G} \cap \on{Gr}_{P_{-}}) \times_{\ol{\on{Gr}}^{\la}_{G} \cap \on{Gr}_{L}} ({\ol{\CW}^{\la}_{\mu}})^{2{\rho}_{G,L}}$ gives us the following data (after the identification of $\on{Gr}_{P_{-},\bar{\mu}}$ and $\on{Gr}_{P,w_{0}({w_{0}}^{L})^{-1}(\bar{\mu})}$ by conjugation via $\tilde{w_{0}}(\tilde{{w_{0}}}^{L})^{-1}$ and using Lemma ~\ref{4} to describe $({\ol{\CW}^{\la}_{\mu}})^{2{\rho}_{G,L}}$):

(a) A $G$-bundle $\CP$ on $\BP^{1}$.

(b) A trivialization $\sig: \CP_{triv}|_{(\BP^{1} \sm 0)} \tilde{\ra} \CP|_{(\BP^{1} \sm 0)\times S}$ such that the corresponding point of $\on{Gr}_{G}(\BC)$ lies in $\ol{\on{Gr}}^{\la}_{G}(\BC)$.

(c) A $L$-bundle $\CP^{L}$ with a trivialization $\sig^{L}: \CP^{L}_{triv}|_{(\BP^{1} \sm 0)} \tilde{\ra} \CP^{L}|_{(\BP^{1} \sm 0)}$.

(d) An isomorphism between $G \times_{L} \CP^{L}$ and $\CP$ such that the corresponding trivializations coincide.

(e) A $P$-bundle $\CP^{P}$ with a trivialization $'\sigma^{P}: \CP^{P}_{triv}|_{(\BP^{1} \sm 0)} \tilde{\ra} \CP^{P}|_{(\BP^{1} \sm 0)}$.

(f) An isomorphism $\psi^{P}$ between $G\times_{P} \CP^{P}$ and $\CP$ such that the corresponding trivializations differ by $\tilde{w_{0}}(\tilde{{w_{0}}}^{L})^{-1}$.

(g) An isomorphism $\psi^{L}$ between $L\times_{P} \CP^{P}$ and $\CP^{L}$.

(h) A $B_{L}$-structure $\phi_{B_{L}}$ on $\CP^{L}$ of degree ${w_{0}}^{L}(\mu)$.

To get a point inside $\ol{\CW}^{\la}_{\mu}$ it is enough to construct a $B$-structure on $\CP^{P}$ of degree $w_{0}(\mu)$.

The $B_{L}$-structure $\phi_{B_{L}}$ on $L\times_{P} \CP^{P}$ induces the desired $B$-structure $\phi_{B}$ on $\CP^{P}$ (by taking the preimage of $\phi_{B_{L}}$ under the natural morphism $\CP^{P} \ra L\times_{P} \CP^{P}$). By the construction the degree of $\phi_{B}$ will be equal to $w_{0}({w_{0}}^{L})^{-1}({w_{0}}^{L}\mu)=w_{0}(\mu)$.

It is easy to see that the point we have constructed lies inside $\ol{R}^{\la}_{G,L,\mu}$.

\epr

Let us consider the following diagram:

\begin{equation}
\label{5}
\xymatrix{\on{Gr}_{G} & \ar[l]^{g^{\la}_{G}} \ol{\on{Gr}}_{G}^{\la} &  \ol{\CW}^{\la}_{\mu} \ar[l]^{\on{\bp}^{\la}_{\mu}}\\
 \on{Gr}_{P_{-},\bar{\mu}}  \ar[u]^{\iota_{\bar{\mu}}} \ar[d]^{\pi_{\bar{\mu}}} & \ar[l] \ol{\on{Gr}}^{\la}_{G} \cap \on{Gr}_{P_{-}, \bar{\mu}} \ar[u]^{\tilde{\iota}^{\la}_{\bar{\mu}}} \ar[d]^{\tilde{\pi}^{\la}_{\bar{\mu}}} & \ol{R}^{\la}_{G,L, \mu} \ar[l]^{\on{\tilde{\bp}}^{\la}_{G,L,\mu}} \ar[u]^{\iota^{\la}_{G,L, \mu}} \ar[d]^{\pi^{\la}_{G,L,\mu}} \\ \on{Gr}_{L,\bar{\mu}} & \ar[l]
 \ol{\on{Gr}}_{L}^{\tilde{\la}} & ({\ol{\CW}^{\la}_{\mu}})^{2{\rho}_{G,L}} \ar[l]^{\on{\tilde{\tilde{\bp}}}^{\la}_{G,L,\mu}} \ar[d]^{\varkappa^{\la}_{\mu,G,L}}\\
&& \ar[ul]^{\on{\bp}^{\tilde{\la}}_{\mu,L}} {\ol{\CW}^{\tilde{\la}}_{\mu,L}}
}
\end{equation}
where $\tilde{\on{\bp}}^{\la}_{G,L,\mu}$ and $\on{\tilde{\tilde{\bp}}}^{\la}_{G,L,\mu}$ are the restrictions of $\on{\bp}^{\la}_{\mu}$ on $\ol{R}^{\la}_{G,L,\mu}$ and $({\ol{\CW}^{\la}_{\mu}})^{2{\rho}_{G,L}}$ respectively, while $\tilde{\iota}^{\la}_{\bar{\mu}}$ and $\tilde{\pi}^{\la}_{\bar{\mu}}$ are the restrictions of $\iota^{\la}_{\bar{\mu}}$ and $\pi^{\la}_{\bar{\mu}}$ to $\ol{\on{Gr}}^{\la}_{G} \cap \on{Gr}_{P_{-}, \bar{\mu}}$.

\lem{Res} \label{14}
$$(\on{\bp}^{\tilde{\la}}_{\mu,L})^{*}[-2\check{\rho}_{L}(\mu)](\tilde{\pi}^{\la}_{\bar{\mu}})_{*}(\tilde{\iota}^{\la}_{\tilde{\mu}})^{!}[2\check{\rho}_{G,L}(\mu)](\on{IC}(\ol{\on{Gr}}^{\la}_{G}))=\bigoplus\limits_{\nu \in \La^{+,\la}_{\mu,G,L}}(\ell^{\nu,\tilde{\la}}_{\mu,G})_{*}\on{IC}(\ol{\CW}^{\nu}_{\mu, L})^{\oplus{n_{\nu}}}.$$
\elem

\prf
By \cite[Proposition 5.3.29]{BD}, the hyperbolic restriction functor $(\pi_{\tilde{\mu}})_{*}(\iota_{\bar{\mu}})^{!}[2\check{\rho}_{G,L}(\mu)]$ sends $(g^{\la}_{G})_{*}\on{IC}(\ol{\on{Gr}^{\la}_{G}})$ to the sum $\bigoplus\limits_{\nu \in \La^{\la,+}_{\bar{\mu},G,L}}{(g^{\nu}_{L})_{*}\on{IC}(\ol{\on{Gr}}^{\nu}_{L})}^{\oplus{n_{\nu}}}$. From the definitions it follows that the perverse sheaf $\bigoplus\limits_{\nu \in \La^{\la,+}_{\bar{\mu},G,L}}{(g^{\nu}_{L})_{*}\on{IC}(\ol{\on{Gr}}^{\nu}_{L})}^{\oplus{n_{\nu}}}$ is supported on $\ol{\on{Gr}}_{L}^{\tilde{\la}}$. So the hyperbolic restriction functor $(\tilde{\pi}^{\la}_{\bar{\mu}})_{*}(\tilde{\iota}^{\la}_{\bar{\mu}})^{!}[2\check{\rho}_{G,L}(\mu)]$ sends $\on{IC}(\ol{\on{Gr}}^{\la}_{G})$ to the sum $\bigoplus\limits_{\nu \in \La^{\la,+}_{\bar{\mu},G,L}}{(h^{\nu, \tilde{\la}}_{L})_{*}\on{IC}(\ol{\on{Gr}}^{\nu}_{L})}^{\oplus{n_{\nu}}}$.

From \cite[Proposition 12.1 c), Proposition 12.4]{FM} and Lemma ~\ref{4} it follows that $(\on{\bp}^{\tilde{\la}}_{\mu,L})^{*}[-2\check{\rho}_{L}(\mu)](\bigoplus\limits_{\nu \in \La^{\la,+}_{\bar{\mu},G,L}}{(h^{\nu, \tilde{\la}}_{L})_{*}\on{IC}(\ol{\on{Gr}}^{\nu}_{L})}^{\oplus{n_{\nu}}})=\bigoplus\limits_{\nu \in \La^{\la,+}_{\bar{\mu},G,L}}(\ell^{\nu,\tilde{\la}}_{\mu,G})_{*}\on{IC}(\ol{\CW}^{\nu}_{\mu, L})^{\oplus{n_{\nu}}}$.
\epr

\ssec{Res}{Proof of Theorem ~\ref{2}}

From Lemma ~\ref{14} and \cite[Proposition 12.1 c), Proposition 12.4]{FM} it follows that to prove Theorem ~\ref{2} it is enough to
constuct the following isomorphism $$(\on{\bp}^{\tilde{\la}}_{\mu,L})^{*}[-2\check{\rho}_{L}(\mu)](\tilde{\pi}^{\la}_{\bar{\mu}})_{*}(\tilde{\iota}^{\la}_{\bar{\mu}})^{!}\on{IC}(\ol{\on{Gr}}^{\la}_{G}))[2\check{\rho}_{G,L}(\mu)]
\stackrel{\sim}{\longrightarrow}$$ $$\stackrel{\sim}{\longrightarrow}
(\varkappa^{\la}_{\mu,G,L})_{*}(\pi^{\la}_{G,L,\mu})_{*}(\iota^{\la}_{G,L,\mu})^{!}(\on{\bp}^{\la}_{\mu})^{*}[-2\check{\rho}_{G}(\mu)]\on{IC}(\ol{\on{Gr}}^{\la}_{G})).$$
Indeed,
$$(\varkappa^{\la}_{\mu,G,L})_{*}(\pi^{\la}_{G,L,\mu})_{*}(\iota^{\la}_{G,L,\mu})^{!}(\on{\bp}^{\la}_{\mu})^{*}\on{IC}(\ol{\on{Gr}}^{\la}_{G}))\is $$
$$(\varkappa^{\la}_{\mu,G,L})_{*}(\pi^{\la}_{G,L,\mu})_{*}(\iota^{\la}_{G,L,\mu})^{!}(\on{\bp}^{\la}_{\mu})^{!}\on{IC}(\ol{\on{Gr}}^{\la}_{G}))[4\check{\rho}_{G}(\mu)]\is $$
$$(\varkappa^{\la}_{\mu,G,L})_{*}(\pi^{\la}_{G,L,\mu})_{*}(\tilde{\on{\bp}}^{\la}_{G,L,\mu})^{!}(\tilde{\iota}^{\la}_{\bar{\mu}})^{!}\on{IC}(\ol{\on{Gr}}^{\la}_{G}))[4\check{\rho}_{G}(\mu)]\is $$
$$(\varkappa^{\la}_{\mu,G,L})_{*}(\tilde{\tilde{\on{\bp}}}^{\la}_{G,L,\mu})^{!}(\tilde{\pi}^{\la}_{\bar{\mu}})_{*}(\tilde{\iota}^{\la}_{\bar{\mu}})^{!}\on{IC}(\ol{\on{Gr}}^{\la}_{G}))[4\check{\rho}_{G}(\mu)]\is $$
$$(\varkappa^{\la}_{\mu,G,L})_{*}(\tilde{\tilde{\on{\bp}}}^{\la}_{G,L,\mu})^{*}(\tilde{\pi}^{\la}_{\bar{\mu}})_{*}(\tilde{\iota}^{\la}_{\bar{\mu}})^{!}\on{IC}(\ol{\on{Gr}}^{\la}_{G}))[4\check{\rho}_{G,L}(\mu)]\is $$
$$(\on{\bp}^{\tilde{\la}}_{\mu,L})^{*}(\tilde{\pi}^{\la}_{\bar{\mu}})_{*}(\tilde{\iota}^{\la}_{\bar{\mu}})^{!}\on{IC}(\ol{\on{Gr}}^{\la}_{G}))[4\check{\rho}_{G,L}(\mu)].$$

The first and fourth isomorphisms follow from the non-characteristic property of the maps $\on{\bp}^{\la}_{\mu}, \on{\bp}^{\tilde{\la}}_{\mu,L}$ (see \cite[Proposition 12.1 c), Proposition 12.4]{FM}) and Lemma ~\ref{4}. The second isomorphism follows from the commutativity of the diagram ~\ref{5}. The third isomorphism follows from Lemma \ref{6} and the base change. The fifth isomorphism follows from Lemma ~\ref{4}.
\qed

\sec{Cs}{Proof of the first part of Theorem ~\ref{7} (2)} \label{32}

Recall that $\bB^{G}(\la):=\coprod\limits_{\mu \in \La_{G}}\bB^{G}(\la)_{\mu}$ where $\bB^{G}(\la)_{\mu}:=\on{Irr}(\ol{R}^{\la}_{\mu})$.
Let us give a geometric construction of a crystal structure on the set $\bB^{G}(\la)$ purely in terms of generalized slices.

Let us take $\la \in \Lap_{G}$, $\mu \in \La_{G}$, $\nu \in \La^{\la,+}_{\mu,G,L}$. Let us denote $\ol{R}^{\la,\nu}_{G,L,\mu}:=(\pi^{\la}_{G,\mu})^{-1}\imath^{\nu}_{\mu,G,L}(\ol{R}^{\nu}_{L,\mu})$.

Let us denote by $\bB^{G}_{L}(\la)_{\nu}$ the set of irreducible components of $\ol{R}^{\la,\nu}_{G,L,\mu}$ of the top dimension.

\prop{CS}
\label{48} For every $\la \in \Lap_{G}$ and $\mu \in \La_{G}$ there is a canonical bijection:

$$\bd^{G}_{L}: \coprod \limits_{\nu \in \La^{\la,+}_{\nu,G,L}} \bB^{G}_{L}(\la)_{\nu} \times \bB^{L}(\nu)_{\mu} \simeq \bB^{G}(\la)_{\mu}$$
\eprop

\rem{CS}
Proposition ~\ref{48} is a reformulation of \cite[Proposition 3.1]{BG} in terms of generalized slices.
\erem

Before proving Proposition ~\ref{48} let us recall some notations of Subsection ~\ref{1}.

Let $P_{-} \subset G$ be a parabolic subgroup of $G$ containing $B_{-}$ and let $U_{P_{-}}$ be the unipotent radical of $P_{-}$, and let $L$ be the corresponding Levi factor. Let us denote by $Z(\check{L})$ the center of $\check{L}$ and by $\La_{G,P_{-}}$ the character lattice of $Z(\check{L})$.Note that $(\on{Gr}_{G})^{Z(L)}=\on{Gr}_{L}$. Recall that $2{\rho}_{G,L}$ is a "subregular" dominant (i.e. the corresponding element of $\La_{G}$ lies inside $\Lap_{G}$ and $\langle 2{\rho}_{G,L}, \al_{i} \rangle > 0$ for any $i \notin I_{L}$) coweight of $Z(L)$. Then $\on{Gr}_{L}$ can be identified with $(\on{Gr}_{G})^{2{\rho}_{G,L}}$.

As in Subsection ~\ref{1} for every $\theta\in\La_{G,P_{-}}$ we have the following diagram (where $\iota_{\theta}$ is a locally closed embedding and $\pi_{\theta}$ is the repelling map of the $\BC^{*}$-action by $2{\rho}_{G,L}$):

\begin{equation}
\xymatrix{ \on{Gr}_{G} & \ar[l]_{\iota_{\theta}} \on{Gr}_{P_{-},\theta} \ar[r]^{\pi_{\theta}} & \on{Gr}_{L,\theta}
}
\end{equation}

Let $T^{\theta}_{P_{-}}:=\iota_{\theta}(\on{Gr}_{P_{-},\theta})$. Take $\la \in \Lap_{G}$, $\nu \in \Lap_{L}$. Let $\theta:=\alpha_{G,P_{-}}(\nu))$. Let $T^{\nu}_{P}:=\iota_{\theta}(\pi_{\theta}^{-1}(\on{Gr}^{\nu}_{L})$. Let $\ol{T}^{\la,\nu}_{G,L,\mu}:=\iota_{\theta}(\pi_{\theta}^{-1}(T^{\nu}_{\mu,L}))\cap \ol{\on{Gr}}^{\la}_{G}$.

Now let us prove Proposition ~\ref{48}.

\prf
We claim that the bijection $\bd^{G}_{L}$ can be uniquely described as follows: one has $\bd^{G}_{L}(\bb_{1},\bb_{2})=\bb$ for $\bb_{1}\in \bB^{G}_{L}(\la)_{\nu}$, $\bb_{2} \in \bB^{L}(\nu)_{\mu}$ if and only if the following conditions hold.

(a) $\bb_{2}$ is a dense subset of $\pi^{\la}_{G,L,\mu}(\bb)$

(b) $(\pi^{\la}_{G,L,\mu})^{-1}(\bb_{2}) \cap \bb_{1}$ is an open dense subset of $\bb$.

To see this note that from Lemma ~\ref{6} and Theorem ~\ref{7} (1) it follows that the map $\bp^{\la}_{G,\mu}$ induces an isomorphism between $\ol{R}^{\la,\nu}_{G,L,\mu}$ and $\ol{T}^{\la,\nu}_{G,L,\mu}$, thus the set $\bB^{G}_{L}(\la)_{\nu}$ can be canonically identified with $\on{Irr}(\ol{T}^{\la,\nu}_{G,L,\mu})$. The set $\on{Irr}(\ol{T}^{\la,\nu}_{G,L,\mu})$ can be canonically identified with the set $\on{Irr}(T^{\nu}_{P} \cap \ol{\on{Gr}}^{\la}_{G})$ (it follows from \cite[3.1]{BG} and \cite{MV}). Now the desired claim directly follows from Theorem ~\ref{7}~(1) and the proof of \cite[Proposition 3.1]{BG}.
\epr

\rem{CS}
Braverman and Gaitsgory consider attractors in their paper while we consider repellents.
\erem

\ssec{op}{Operators $\on{e}_{i}$ and $\on{f}_{i}$} \label{31}  Let $P_{i}$ be the sub-minimal parabolic subgroup containing $B_{-}$ which corresponds to the simple root $\check{\alpha}_{i}$. Let $L_{i}$ be the corresponding Levi factor. Let us consider the bijection of Proposition ~\ref{48} for $L=L_{i}$. Note that for every $\mu \in \La_{G}$, $\nu \in \Lap_{L_{i}}$, the multiplicity of $\mu$ in the irreducible representation of $\check{L}_{i}$ of highest weight $\nu$ is 0 or 1 (it follows from the representation theory of $SL_{2}(\BC)$). So the set $\bB^{L_{i}}(\nu)_{\mu}$ contains no more than one element. Take $\bb \in \bB^{G}(\la)_{\mu}$. Let us assume $\bb=\bd^{G}_{L_{i}}(\bb_{1}, \bb_{2})$, $\bb_{1} \in \bB^{G}_{L_{i}}(\la)_{\nu}$, $\bb_{2} \in \bB^{L_{i}}(\nu)_{\mu}$. Operations $\on{e}_{i}$ and $\on{f}_{i}$ are defined as follows:

 \begin{equation} \on{e}_{i} \cdot \bb =
\begin{cases}
\bd^{G}_{L_{i}}(\bb_{1}, \hat{\bb}_{2}), \hspace{0,1cm} \on{if} \on{there} \on{exists} \hat{\bb}_{2} \in \bB^{G}_{L_{i}}(\la)_{\nu+\alpha_{i}}\\
0 \hspace{0,1cm} \on{otherwise}.
\end{cases}
\end{equation}

 \begin{equation} \on{f}_{i} \cdot \bb =
\begin{cases}
\bd^{G}_{M_{i}}(\bb_{1}, \hat{\bb}_{2}), \hspace{0,1cm} \on{if} \on{there} \on{exists} \hat{\bb}_{2} \in \bB^{G}_{L_{i}}(\la)_{\nu-\alpha_{i}}\\
0 \hspace{0,1cm} \on{otherwise}.
\end{cases}
\end{equation}

Now the first part of Theorem ~\ref{7} (2) follows from \cite[Subsection 3.3]{BG} and \cite[Theorem 3.1 (1), (2)]{BG}.

\sec{L2}{Proof of Theorem ~\ref{7} (2)} \label{23}

\ssec{Mtr}{Matrix description of the multiplication.}
\prop{Mtr} \cite[5.9]{FKPRW}
The multiplication morphism $\ol{\CW}^{\la_{1}}_{\mu_{1}} \times \ol{\CW}^{\la_{2}}_{\mu_{2}} \ra \ol{\CW}^{\la_{1}+\la_{2}}_{\mu_{1}+\mu_{2}}$ is given by $(g_{1}, g_{2}) \mapsto \pi(g_{1} g_{2})$, $g_{1} \in \ol{\CW}^{\la_{1}}_{\mu_{1}}, g_{2} \in \ol{\CW}^{\la_{2}}_{\mu_{2}}$ where

$ \pi : U[z]$ $\backslash$  $U((z^{-1}))T[[z^{-1}]]_{1}z^{\mu}U((z^{-1}))/U[z] \simeq U_{1}[[z^{-1}]]T_{1}[[z^{-1}]]z^{\mu}U_{1}[[z^{-1}]].$
\eprop

\ssec{Split}{Splitting of convolution}

Let $T_{\mu_{1}}\s T_{\mu_{2}}:=U_{-}(\CK)z^{\mu_{1}}B_{-}(\CO) \times_{B_{-}(\CO)} T_{\mu_{2}}$.

Let $T^{\la_{1}}_{\mu_{1}}\s T^{\la_{2}}_{\mu_{2}}:=(T_{\mu_{1}}\s T_{\mu_{2}})\cap (\on{Gr}^{\la_{1}}_{G}\s \on{Gr}^{\la_{2}}_{G}).$

Let $\ol{T}^{\la_{1}}_{\mu_{1}}\s \ol{T}^{\la_{2}}_{\mu_{2}}:=(T^{\la_{1}}_{\mu_{1}}\s T^{\la_{2}}_{\mu_{2}})\cap (\ol{\on{Gr}}^{\la_{1}}_{G}\s \ol{\on{Gr}}^{\la_{2}}_{G}).$

\lem{L2} \label{13} There exists an isomorphism $\tau^{\la_{1},\la_{2}}_{\mu_{1},\mu_{2}}$ between $T^{\la_{1}}_{\mu_{1}} \times T^{\la_{2}}_{\mu_{2}}$ and $T^{\la_{1}}_{\mu_{1}}\s T^{\la_{2}}_{\mu_{2}}$ which induces an identity morphism on the set of irreducible components of maximal dimension. Here we identify $\on{Irr}(T^{\la_{1}}_{\mu_{1}} \s T^{\la_{2}}_{\mu_{2}}) \simeq \bB^{G}(\la_{1})_{\mu_{1}} \times \bB^{G}(\la_{2})_{\mu_{2}} \simeq \on{Irr}(T^{\la_{1}}_{\mu_{1}} \times T^{\la_{2}}_{\mu_{2}})$.
\elem

\prf One has a natural isomorphism $\on{\on{t}}_{\mu}:T_{\mu}\simeq U[z^{-1}]_{1}$. To see this let us take a point $[uz^{\mu}] \in T_{\mu}$ and rewrite it as $z^{\mu}\tilde{u}$ and then apply $\xi^{-1}$ to $\tilde{u}$.

Let us start with constructing an isomorphism $\tau_{\mu_{1}, \mu_{2}}$ between $T_{\mu_{1}} \times T_{\mu_{2}}$ and $T_{\mu_{1}} \s T_{\mu_{2}}$. Recall that $T_{\mu_{1}} \s T_{\mu_{2}} = U_{-}(\CK)z^{\mu_{1}}B_{-}(\CO) \times_{B_{-}(\CO)} T_{\mu_{2}} =  U_{-}(\CK)z^{\mu_{1}} U_{-}(\CO) \times_{U_{-}(\CO)} T_{\mu_{2}}$.
For a point $([u_{1}z^{\mu_{1}}], [u_{2}z^{\mu_{2}}]) \in T_{\mu_{1}} \times T_{\mu_{2}}$ let us denote $u^{-}_{1} := \on{t}_{\mu_{1}}([u_{1}z^{\mu_{1}}]).$ Directly from the definition of $\on{\on{t}}_{\mu}$ we see that $([u_{1}z^{\mu_{1}}], [u_{2}z^{\mu_{2}}]) = ([z^{\mu_{1}}u^{-}_{1}], [u_{2}z^{\mu_{2}}])$. Let us define $\tau_{\mu_{1}, \mu_{2}}([u_{1}z^{\mu_{1}}], [u_{2}z^{\mu_{2}}]) := (z^{\mu_{1}}u^{-}_{1}, [u_{2}z^{\mu_{2}}])$.

Let us construct the inverse morphism $\tau_{\mu_{1},\mu_{2}}^{-1}:T_{\mu_{1}}\s T_{\mu_{2}}\ra T_{\mu_{1}} \times T_{\mu_{2}}$.$$T_{\mu_{1}}\s T_{\mu_{2}}=U_{-}(\CK)z^{\mu_{1}} U_{-}(\CO) \times_{U_{-}(\CO)} T_{\mu_{2}}=z^{\mu_{1}}U_{-}(\CK)\times_{U_{-}(\CO)} T_{\mu_{2}}.$$

Let us take a point $(z^{\mu_{1}}u_{1},[u_{2}z^{\mu_{2}}])\in T_{\mu_{1}}\s T_{\mu_{2}}$. Let $\on{\bf{mult}}^{-1}(u_{1})=:(u_{1}^{-},u_{\CO})$. By the definition $u_{1}^{-}\in U_{-}[z^{-1}]_{1}$, $u_{\CO} \in U_{-}(\CO)$. Let $\tau_{\mu_{1},\mu_{2}}^{-1}(z^{\mu_{1}}u_{1},[u_{2}z^{\mu_{2}}]):=([z^{\mu_{1}}u_{1}^{-}],[u_{\CO}u_{2}z^{\mu_{2}}]).$ It is easy to see that $\tau_{\mu_{1},\mu_{2}}$ and $\tau_{\mu_{1},\mu_{2}}^{-1}$ are mutually inverse.

It follows from the definitions that $\tau_{\mu_{1}, \mu_{2}}$ restricts to the isomorphism $$\tau^{\la_{1}, \la_{2}}_{\mu_{1}, \mu_{2}}:  T^{\la_{1}}_{\mu_{1}} \times T^{\la_{2}}_{\mu_{2}} \simeq T^{\la_{1}}_{\mu_{1}}\s T^{\la_{2}}_{\mu_{2}}.$$

Before proving the last statement of the lemma let us first recall the isomorphism (let us denote it $\Upsilon^{\la_{1},\la_{2}}_{\mu_{1},\mu_{2}}$) between $\bB^{G}(\la_{1})_{\mu_{1}} \times \bB^{G}(\la_{2})_{\mu_{2}}$ and $\on{Irr}(T^{\la_{1}}_{\mu_{1}}\s T^{\la_{2}}_{\mu_{2}})$. Let $p^{\la_{1},\la_{2}}_{\mu_{1},\mu_{2}}:T^{\la_{1}}_{\mu_{1}}\s T^{\la_{2}}_{\mu_{2}}\ra T^{\la_{1}}_{\mu_{1}}$ denote the projection onto the first factor. The isomorphism $\Upsilon^{\la_{1},\la_{2}}_{\mu_{1},\mu_{2}}$ can be uniquely characterized as follows: one has $\Upsilon^{\la_{1},\la_{2}}_{\mu_{1},\mu_{2}}(\on{\bb_{1}},\on{\bb_{2}})=\on{\bb}$ for $\on{\bb_{1}} \in \bB^{G}(\la_{1})_{\mu_{1}}$, $\on{\bb_{2}} \in \bB^{G}(\la_{2})_{\mu_{2}}$ if and only if the following conditions hold.

1) $\on{\bb_{1}}$ is a dense subset of $p^{\la_{1},\la_{2}}_{\mu_{1},\mu_{2}}(\on{\bb}).$

2) $(p^{\la_{1},\la_{2}}_{\mu_{1},\mu_{2}})^{-1}(\on{\bb_{1}})\cap \on{\bb_{2}}$ is a dense subset of $\on{\bb}.$

Let $pr^{\la_{1},\la_{2}}_{\mu_{1},\mu_{2}}$ denote the projection of $T^{\la_{1}}_{\mu_{1}} \times T^{\la_{2}}_{\mu_{2}}$ to the first factor. Directly from the definitions
we see that the following diagram is commutative:

\begin{equation}
\xymatrix{ T^{\la_{1}}_{\mu_{1}} \times T^{\la_{2}}_{\mu_{2}}  \ar[rr]^{\tau^{\la_{1},\la_{2}}_{\mu_{1},\mu_{2}}} \ar[dr]_{pr^{\la_{1},\la_{2}}_{\mu_{1},\mu_{2}}} &&  T^{\la_{1}}_{\mu_{1}} \s T^{\la_{2}}_{\mu_{2}} \ar[ld]^{p^{\la_{1},\la_{2}}_{\mu_{1},\mu_{2}}} \\
& T^{\la_{1}}_{\mu_{1}}
}
\end{equation}

So to prove the last statement of the дemma it is enough to show that for any point $x \in  T^{\la_{1}}_{\mu_{1}}$ the map $\tau^{\la_{1},\la_{2}}_{\mu_{1},\mu_{2}}$ induces the identity morphism between $\on{Irr}({pr^{\la_{1},\la_{2}}_{\mu_{1},\mu_{2}}}^{-1}(x))$ and $\on{Irr}({p^{\la_{1},\la_{2}}_{\mu_{1},\mu_{2}}}^{-1}(x))$ (here we identify both of these sets with $\bB(\la_{2})_{\mu_{2}}$).

Let $x=[z^{\mu_{1}}u^{-}_{1}]$, $u^{-}_{1} \in U_{-}[z^{-1}]_{1}$. Then ${pr^{\la_{1},\la_{2}}_{\mu_{1},\mu_{2}}}^{-1}(x)=[z^{\mu_{1}}u^{-}_{1}] \times T^{\la_{2}}_{\mu_{2}}$ and $\tau^{\la_{1},\la_{2}}_{\mu_{1},\mu_{2}}$ sends $([z^{\mu_{1}}u^{-}_{1}],y) \in [z^{\mu_{1}}u^{-}_{1}] \times T^{\la_{2}}_{\mu_{2}}$ to $(z^{\mu_{1}}u^{-}_{1},y) \in {p^{\la_{1},\la_{2}}_{\mu_{1},\mu_{2}}}^{-1}(x)$. Thus the statement follows.

\epr

Let us recall the geometric construction of the crystal $\on{\bB^{G}}(\la_{1}) \otimes \on{\bB^{G}}(\la_{2})$ from \cite[Subsection 5.1]{BG}:

\prop{CS}
There exists an isomorphism of crystals $$\Xi^{\la_{1},\la_{2}}:\on{\bB^{G}}(\la_{1}) \otimes \on{\bB^{G}}(\la_{2}) \simeq \coprod \limits_{\la_{3} \in \Lap_{G}} \on{Irr}(\on{Gr}^{\la_{1}}\s \on{Gr}^{\la_{2}})^{\la_{3}} \times \bB^{G}(\la_{3})$$
where the crystal structure in the RHS is induced from the crystal structure of $\bB^{G}(\la_{3})$.
\eprop

The isomorphism $\Xi^{\la_{1},\la_{2}}$ on the level of sets can be constructed as follows.
Note that $\on{\bB}^{G}(\la_{1}) \otimes \on{\bB}^{G}(\la_{2}) \simeq \coprod \limits_{\mu_{1}, \mu_{2} \in \La_{G}} \on{Irr}(T^{\la_{1}}_{\mu_{1}}\s T^{\la_{2}}_{\mu_{2}})$ and $\coprod \limits_{\la_{3} \in \Lap_{G}} \on{Irr}(\on{Gr}^{\la_{1}}\s \on{Gr}^{\la_{2}})^{\la_{3}} \times \bB^{G}(\la_{3}) \simeq \coprod \limits_{\mu \in \La_{G}} \on{Irr}(\on{m}^{-1}(T_{\mu}) \cap (\on{Gr}^{\la_{1}}\s \on{Gr}^{\la_{2}})).$

Let us take $\on{\bb} \in \on{Irr}(T^{\la_{1}}_{\mu_{1}}\s T^{\la_{2}}_{\mu_{2}})$. The morphism $\Xi^{\la_{1},\la_{2}}$ sends $\on{\bb}$ to its closure in the variety  $\on{\bf{m}}^{-1}(T_{\mu_{1}+\mu_{2}}) \cap (\on{Gr}^{\la_{1}}\s \on{Gr}^{\la_{2}}).$

Let $\on{\bf{m}}^{\la_{1}, \la_{2}}_{\mu_{1},\mu_{2}}: T^{\la_{1}}_{\mu_{1}}\s T^{\la_{2}}_{\mu_{2}} \ra \ol{T}^{\la_{1}+\la_{2}}_{\mu_{1}+\mu_{2}}$ and $\ol{\on{\bf{m}}}^{\la_{1}, \la_{2}}_{\mu_{1},\mu_{2}}: \ol{T}^{\la_{1}}_{\mu_{1}}\s \ol{T}^{\la_{2}}_{\mu_{2}} \ra \ol{T}^{\la_{1}+\la_{2}}_{\mu_{1}+\mu_{2}}$ be the convolution morphisms.

\lem{L2} \label{12} The map $\on{\bf{m}}^{\la_{1}, \la_{2}}_{\mu_{1},\mu_{2}}$ induces the retraction of crystals, i.e. under the identification of $\coprod \limits_{\mu_{1}, \mu_{2} \in \La_{G}}\on{Irr}(T^{\la_{1}}_{\mu_{1}}\s T^{\la_{2}}_{\mu_{2}}) \is \on{\bB}^{G}(\la_{1}) \otimes \on{\bB}^{G}(\la_{2})$
with $\coprod \limits_{\la_{3} \in \Lap_{G}}\on{Irr}(\on{Gr}^{\la_{1}}*\on{Gr}^{\la_{2}})^{\la_{3}} \times \bB^{G}(\la_{3})$ via $\Xi^{\la_{1},\la_{2}}$, $\on{\bf{m}}^{\la_{1}, \la_{2}}_{\mu_{1},\mu_{2}}$ maps $\on{\bb} \in \on{Irr}(T^{\la_{1}}_{\mu_{1}}\s T^{\la_{2}}_{\mu_{2}})$ to itself if $\on{\bb} \in \bB^{G}(\la_{1}+\la_{2})$, and to zero otherwise.
\elem

\prf Let us take an irreducible component $\on{\bb} \in \on{Irr}(T^{\la_{1}}_{\mu_{1}}\s T^{\la_{2}}_{\mu_{2}})$. Let us suppose $\on{\bb} \notin \bB(\la_{1}+\la_{2})_{\mu_{1}+\mu_{2}}$. Then there exists $\la_{3} < \la_{1}+\la_{2}$ such that $\on{\bb} \subset {\on{\bf{m}}^{\la_{1}, \la_{2}}_{\mu_{1},\mu_{2}}}^{-1}(T_{\mu_{1}+\mu_{2}}) \cap (\on{Gr}^{\la_{1}}_{G}\s \on{Gr}^{\la_{2}}_{G})^{\la_{3}}$, thus $\on{\bf{m}}^{\la_{1}, \la_{2}}_{\mu_{1},\mu_{2}}(\on{b}) \subset \ol{T}^{\la_{3}}_{\mu_{1}+\mu_{2}}$, so $\on{dim}\on{\bf{m}}^{\la_{1}, \la_{2}}_{\mu_{1},\mu_{2}}(\on{b}) \leq \on{dim}(\ol{T}^{\la_{3}}_{\mu_{1}+\mu_{2}})=\check{\rho}_{G}(\la_{3}-\mu_{1}-\mu_{2})<\check{\rho}_{G}(\la_{1}+\la_{2}-\mu_{1}-\mu_{2})=\on{dim}(\ol{T}^{\la_{1}+\la_{2}}_{\mu_{1}+\mu_{2}})$. Thus $\on{\bf{m}}^{\la_{1}, \la_{2}}_{\mu_{1},\mu_{2}}(\on{\bb})$ is not of maximal dimension, so $\on{\bf{m}}^{\la_{1}, \la_{2}}_{\mu_{1},\mu_{2}}$ maps $\on{\bb}$ to $0$.

Let us suppose $\on{\bb} \in \bB^{G}(\la_{1}+\la_{2})_{\mu_{1}+\mu_{2}}$ i.e. the closure of $\on{\bb}$ inside $\on{\bf{m}}^{-1}(T_{\mu_{1}+\mu_{2}}) \cap (\on{Gr}^{\la_{1}}_{G}\s \on{Gr}^{\la_{2}}_{G})$ actually lies in $\on{\bf{m}}^{-1}(T_{\mu_{1}+\mu_{2}}) \cap \on{\bf{m}}^{-1}(\on{Gr}^{\la_{1}+\la_{2}}_{G})$. From stratified semismallness of the morphism $\on{\bf{m}}$ (see \cite[Lemma 4.4]{MV}) it follows that $\on{dim}(\on{\bf{m}}^{\la_{1}, \la_{2}}_{\mu_{1},\mu_{2}}(\on{\bb}))=\on{dim}(\on{\bb})=\check{\rho}_{G}(\la_{1}+\la_{2}-\mu_{1}-\mu_{2})$ (actually, the map $\on{\bf{m}}$ induces
an isomorphism between $(\on{Gr}^{\la_{1}}_{G}\s \on{Gr}^{\la_{2}}_{G})^{\la_{1}+\la_{2}}=(\on{Gr}^{\la_{1}}_{G}\s \on{Gr}^{\la_{2}}_{G}) \cap{\on{\bf{m}}}^{-1}(\on{Gr}^{\la_{1}+\la_{2}}_{G})$ and $\on{Gr}^{\la_{1}+\la_{2}}_{G})$. Thus the closure of $\on{\bf{m}}^{\la_{1}, \la_{2}}_{\mu_{1},\mu_{2}}(\on{\bb})$ inside $\ol{T}^{\la_{1}+\la_{2}}_{\mu_{1}+\mu_{2}}$ is an element of $\bB^{G}(\la_{1}+\la_{2})_{\mu_{1}+\mu_{2}}$. Directly from the definitions (see \cite[Subsection 5.1]{BG}) we see that $\on{\bf{m}}^{\la_{1}, \la_{2}}_{\mu_{1},\mu_{2}}$ maps $\bb$ to itself.

\epr

\lem{L2} \label{11} The following diagram is commutative:
\begin{equation}
\xymatrix{
\ol{R}^{\la_{1}}_{\mu_{1}} \times \ol{R}^{\la_{2}}_{\mu_{2}} \ar[rr]^{\on{p}^{\la_{1}}_{\mu_{1}}\times\on{p}^{\la_{2}}_{\mu_{2}}}  \ar[d]^{\ol{\bf{\kappa}}^{\la_{1}, \la_{2}}_{\mu_{1}, \mu_{2}}} && \ol{T}^{\la_{1}}_{\mu_{1}} \times \ol{T}^{\la_{2}}_{\mu_{2}} \ar[rr]^{\tau^{\la_{1}, \la_{2}}_{\mu_{1}, \mu_{2}}} &&  \ol{T}^{\la_{1}}_{\mu_{1}}\s \ol{T}^{\la_{2}}_{\mu_{2}} \ar[d]^{\ol{\on{\bf{m}}}^{\la_{1}, \la_{2}}_{\mu_{1}, \mu_{2}}} \\
\ol{R}^{\la_{1}+\la_{2}}_{\mu_{1}+\mu_{2}} \ar[rrrr]^{\on{p}^{\la_{1}+\la_{2}}_{\mu_{1}+\mu_{2}}} &&&& \ol{T}^{\la_{1}+\la_{2}}_{\mu_{1}+\mu_{2}}
}
\end{equation}
\elem

\prf Let us take a point $(z^{\mu_{1}}u_{1}^{-}, z^{\mu_{2}}u_{2}^{-}) \in \ol{R}^{\la_{1}}_{\mu_{1}} \times \ol{R}^{\la_{2}}_{\mu_{2}}$. Then by the definitions it is enough to prove the commutativity of the following diagram:

\begin{equation}
\xymatrix{
(z^{\mu_{1}}u_{1}^{-}, z^{\mu_{2}}u_{2}^{-}) \ar[rr]^{\on{p}^{\la_{1}}_{\mu_{1}}\times\on{p}^{\la_{2}}_{\mu_{2}}}  \ar[d]^{\ol{\bf{\kappa}}^{\la_{1}, \la_{2}}_{\mu_{1}, \mu_{2}}} && ([z^{\mu_{1}}u_{1}^{-}], [z^{\mu_{2}}u_{2}^{-}]) \ar[rr]^{\tau^{\la_{1}, \la_{2}}_{\mu_{1}, \mu_{2}}} &&  (z^{\mu_{1}}u_{1}^{-}, [z^{\mu_{2}}u_{2}^{-}]) \ar[d]^{\ol{\on{\bf{m}}}^{\la_{1}, \la_{2}}_{\mu_{1}, \mu_{2}}} \\
\pi(z^{\mu_{1}}u_{1}^{-}z^{\mu_{2}}u_{2}^{-}) \ar[rrrr]^{\on{p}^{\la_{1}+\la_{2}}_{\mu_{1}+\mu_{2}}} &&&& [z^{\mu_{1}}u_{1}^{-}z^{\mu_{2}}u_{2}^{-}]
}
\end{equation}

In other words, we have to prove that $\on{p}^{\la_{1}+\la_{2}}_{\mu_{1}+\mu_{2}}(\pi(z^{\mu_{1}}u_{1}^{-}z^{\mu_{2}}u_{2}^{-}))=[z^{\mu_{1}}u_{1}^{-}z^{\mu_{2}}u_{2}^{-}]$. It follows from the fact that $\on{p}^{\la_{1}+\la_{2}}_{\mu_{1}+\mu_{2}}(\pi(z^{\mu_{1}}u_{1}^{-}z^{\mu_{2}}u_{2}^{-}))=[\pi(z^{\mu_{1}}u_{1}^{-}z^{\mu_{2}}u_{2}^{-})]=[z^{\mu_{1}}u_{1}^{-}z^{\mu_{2}}u_{2}^{-}]$. The last equality follows from the definition of $\pi$.
\epr

Now Theorem ~\ref{7} (2) follows from Lemmas ~\ref{13}, ~\ref{12}, ~\ref{11}.

\end{document}